\documentclass{elsart} 
\usepackage{graphicx,amssymb,lineno}
\usepackage{natbib}
\newcommand{\rtbpe}{Restricted Three Body Problem }
\newcommand{\rtbp}{Restricted Three Body Problem}
\newcommand{\p}{p_{\Phi}}
\newcommand{\q}{q_{\Phi}}
\begin{document}

\begin{frontmatter}
\title{The role of the unstable equilibrium points in the transfer
of matter in galactic potentials}

\author[label1]{M. Romero-G\'omez}
\ead{merce.romerogomez@oamp.fr}
\author[label2]{J.J. Masdemont}
\author[label3]{C. Garc\'{\i}a-G\'omez}
\author[label1]{E. Athanassoula}

\address[label1]{Laboratoire d'Astrophysique de Marseille (LAM), UMR6110, 
Pole de l'Etoile Site de Chateau-Gombert, 38 rue Frederic Joliot-Curie, 13388 
Marseille C\'edex 13, France}
\address[label2]{I.E.E.C \& Dep. Mat. Aplicada I, Universitat Polit\`ecnica de
Catalunya, Diagonal 647, 08028 Barcelona, Spain}
\address[label3]{D.E.I.M., Universitat Rovira i Virgili, Campus Sescelades, 
Avd. dels Pa\"{\i}sos Catalans 26, 43007 Tarragona, Spain}

\begin{abstract}
We study the role of the unstable equilibrium points in the transfer of matter in a galaxy using 
the potential of a rotating triaxial system. In particular, we study the neighbourhood of these 
points for energy levels and for main model parameters where the zero velocity curves just open 
and form a bottleneck in the region. For these energies, the transfer of matter from the inner 
to the outer parts and vice versa starts being possible. We study how the dynamics around 
the unstable equilibrium points is driven, by performing a partial normal form scheme and by computing 
the invariant manifolds of periodic orbits and quasi-periodic orbits using the reduced Hamiltonian. 
In particular, we compute some homoclinic and heteroclinic orbits playing a crucial role. Our results 
also show that in slow rotating and/or axisymmetric systems the hyperbolic character of the 
equilibrium points is cancelled, so that no transfer of matter is possible through the bottleneck.
\end{abstract}

\end{frontmatter}

\section{Introduction}
\label{sec:intro}
This paper focuses on the study of the dynamics around the hyperbolic equilibrium points, $L_1$ and 
$L_2$, of a non-axisymmetric galactic potential, i.e. a potential whose principal axes on the 
$(x,y)$ plane are different. In particular, we analyse the role of the invariant manifolds 
associated with the unstable periodic orbits and with invariant tori around $L_1$ and $L_2$ in the 
large scale transfer of matter within the system. For this purpose, we use both semi-analytical and 
numerical techniques to compute the invariant manifolds and to study the role they play in the 
global morphology. This is a modern approach using dynamical systems which has successfully been 
applied in celestial mechanics and astrodynamics (e.g. \cite{gom91,jor99,gom01,koo00,gom04}). A similar 
technique has already been used to explain the spiral arm and ring morphology in barred galaxies 
\cite{rom06,rom07,ath08}. Previous theories believe that spiral arms are density waves in a disc 
galaxy \cite{lin63}. The density waves propagate from the centre towards the principal resonances of the 
galaxy, where they damp down \cite{too69}. Other replenishment theories have been proposed, therefore, to obtain 
long-lived spirals (see \cite{ath84} for a review). Using this innovative approach, we obtain 
outer rings and long-lived spirals. This paper is intended to explain the geometrical behaviour of 
such structures in a general manner. 

Such studies have many applications in galactic dynamics. Elliptical galaxies can be triaxial, (i.e. 
their principal axes can be all different from each other) while disc galaxies can contain 
several triaxial components, rotating or non-rotating, and with widely varying degree of 
non-axisymmetry, 
such as haloes, bulges, bars, or oval discs. In this paper, we consider a galactic model that describes 
a triaxial system. By studying its dynamics within a wide range of parameters, we will explain the 
dynamics of all the components and make links with their morphology. The logarithmic potential is a 
very suitable model for our purposes, because it has a simple expression, it is analytic, its series 
expansion can be calculated up to a high order, and it admits both semi-analytical and numerical 
treatments. 

\begin{figure}
\begin{center}
\includegraphics[scale=0.43,angle=-90.0]{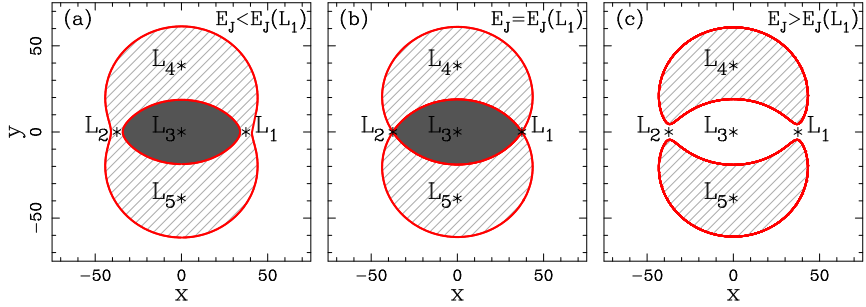}
\caption{Location of the equilibrium points and the zero velocity curves. {\bf (a)} Zero
velocity curves at an energy level smaller than the energy of $L_1$ and $L_2$. Three regions are thus 
defined, namely an inner (solid dark grey), an outer (solid white) and the forbidden region (hatched 
light grey), and no transfer of matter between the inner and outer regions is possible. {\bf (b)} Zero 
velocity curves at the energy level of the unstable equilibrium points $L_1$ and $L_2$. {\bf (c)} Zero 
velocity curves at an energy level larger than the energy of $L_1$ and $L_2$. A bottleneck appears and 
matter can transit from the inner to the outer regions (solid white) and vice versa.}
\label{fig:bottleneck}
\end{center}
\end{figure}

We choose a reference frame such that the origin of coordinates coincides with the centre of the 
galaxy and the triaxial system is fixed, i.e. a reference frame rotating with the system. The idea 
of the present analysis is based on the fact that if we consider a range of energies for which the 
particles are confined in the inner region defined by the zero velocity curve, in our case, for 
energies lower or equal to that of the unstable equilibrium points ($E_J(L_1)$), then there is no 
possible transfer of matter from the inner region of the galaxy to the outer region, or vice versa 
(see Fig.~\ref{fig:bottleneck}(a) and (b)). In this paper, we focus on energies slightly larger than 
that of the unstable equilibrium point, for which an opening in form of a bottleneck appears. We 
will hereafter refer to this aperture as bottleneck. Transfer of matter can be possible through this 
bottleneck (see Fig.~\ref{fig:bottleneck}(c)), and we study the objects that drive the motion in 
this region.

Therefore, our first goal is to study the neighbourhood of the unstable equilibrium points in
the energy range where the bottleneck appears. This energy range is suitable to apply the
normal form technique, which is known as the reduction to the centre manifold. As previously mentioned,
it has been successfully applied to celestial mechanics problems, in particular to the \rtbp, e.g. 
\cite{gom91,gom01,jor99}. The reduced Hamiltonian obtained from the reduction to the centre 
manifold process gives a good qualitative description of the phase space near the equilibrium 
points and uncouples the centre manifold from its hyperbolic behaviour. The procedure is similar 
to the Birkhoff normal form, usually used to study the stability properties around the central 
equilibrium point of the galaxy and to compute the families of periodic orbits around it, e.g. 
\cite{gus66,mir89,bel07}.

Our second goal is to study the behaviour of the hyperbolic invariant manifolds associated 
with the orbits contained in the centre manifold; these are essentially invariant manifolds 
of periodic orbits and invariant tori. As is well-known, stable and unstable invariant 
manifolds are dynamical features, that are responsible for the global dynamics in 
a dynamical system. In order to obtain the global picture of the transfer of mass, it is 
necessary to compute the invariant manifolds and to look for intersections between 
different parts, that is, obtaining the possible heteroclinic and homoclinic connections. These
type of connections have been recently studied in the \rtbpe and applied to obtaining transit
and non-transit orbits in celestial mechanics \cite{koo00,gom04}.

In Section~\ref{sec:mod}, we describe the characteristics of the galactic 
model. In Section~\ref{sec:mot}, we present the equations of motion and we study
the effect of the main parameters on the linear behaviour of the equilibrium points. 
In Section~\ref{sec:red}, we explain in detail the reduction to the centre manifold 
in the particular case of the logarithmic potential and perform a study of the
practical convergence of the reduced Hamiltonian. In Section~\ref{sec:res}, we compute 
the invariant manifolds associated with periodic orbits and quasi-periodic orbits and we study 
the role they play in the transfer of matter. We also perform an analysis on the variation of the 
free parameters of the potential. Finally, in Section~\ref{sec:conc} we conclude.

\section{Galactic model}
\label{sec:mod}

\begin{figure}
\begin{center}
\includegraphics[scale=0.25,angle=-90.0]{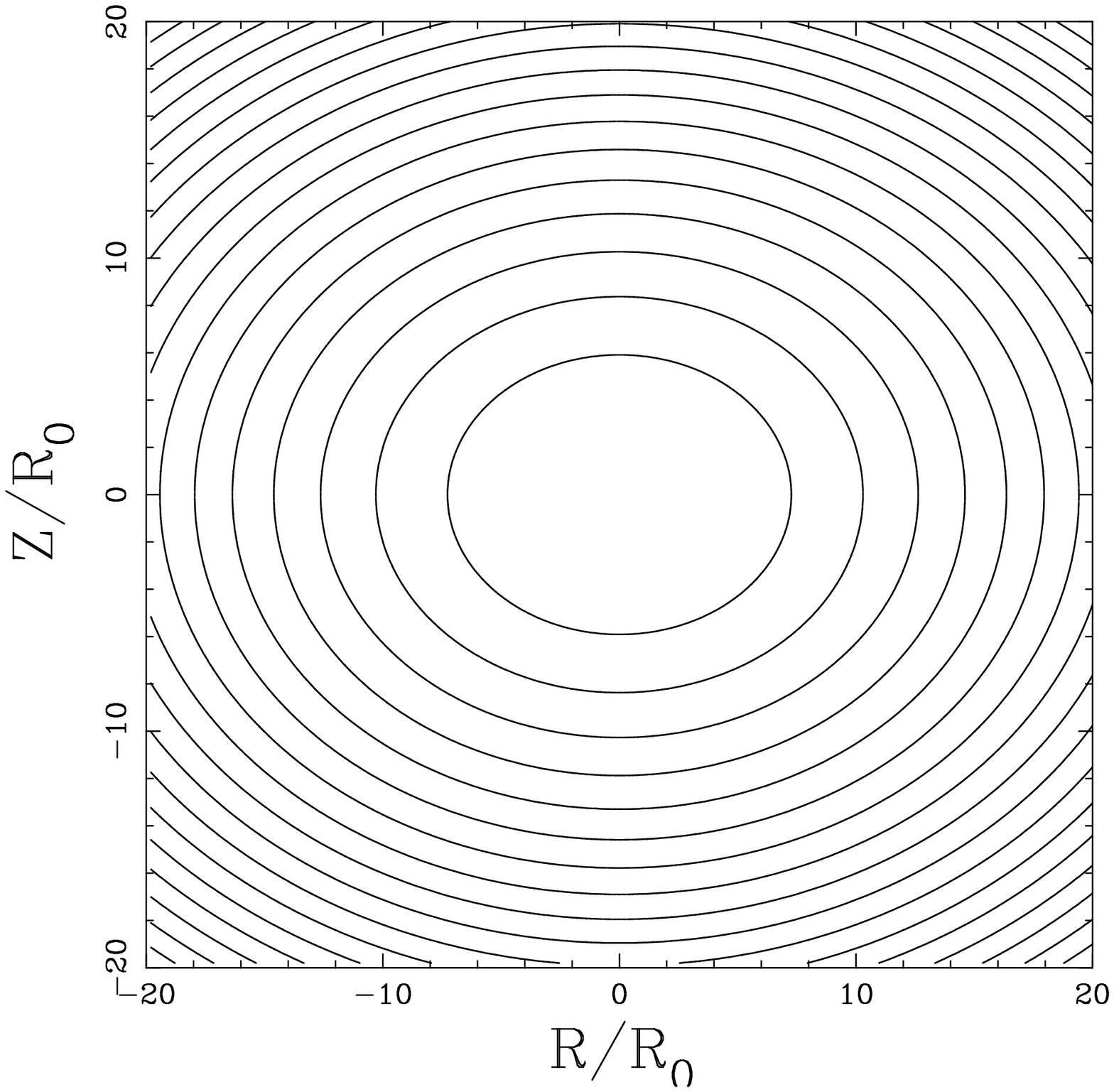}\hspace{0.5cm}
\includegraphics[scale=0.25,angle=-90.0]{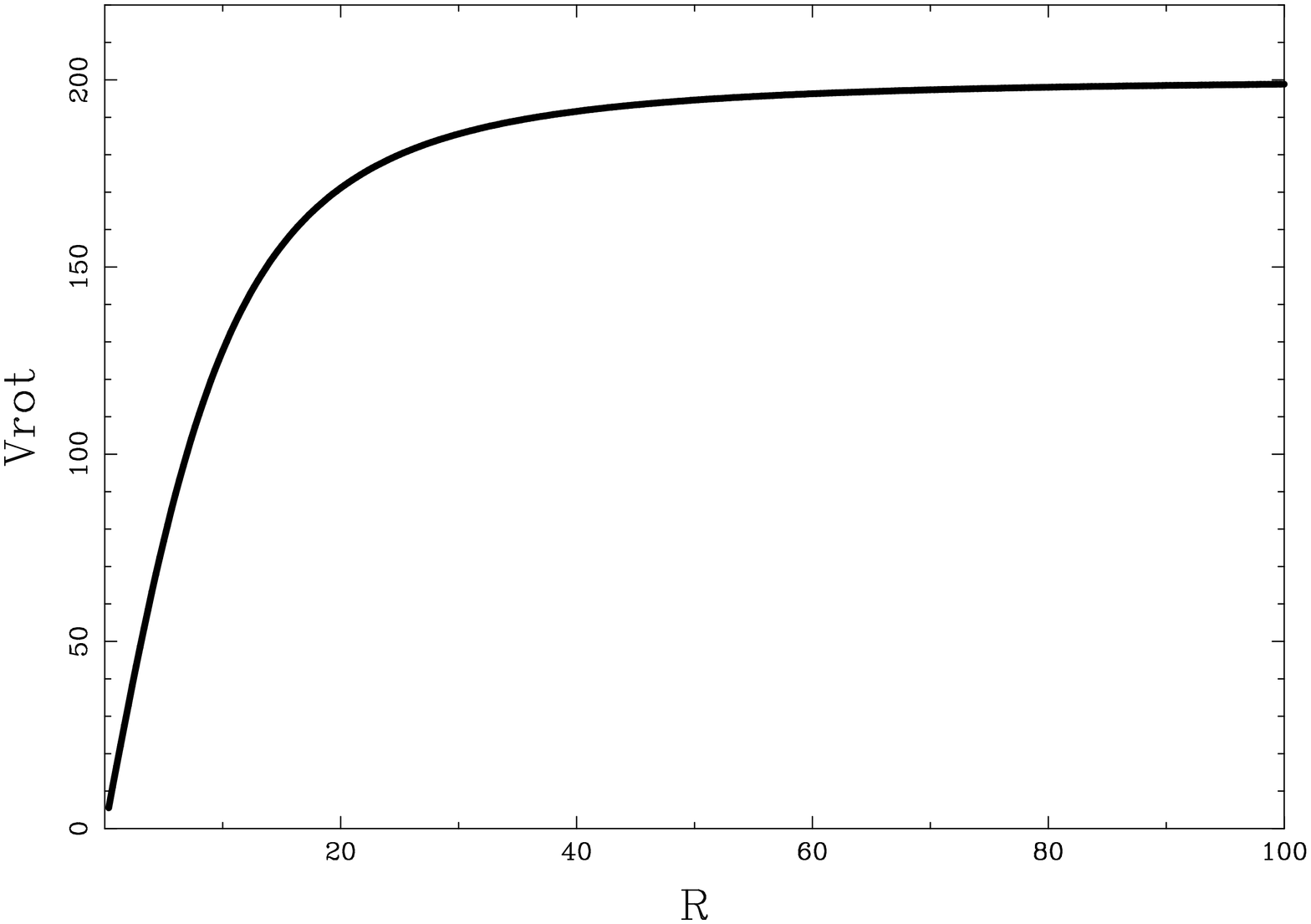}
\caption{Characteristics of the logarithmic potential. {\sl Left panel:} Isodensity curves 
for ``model 1'' with parameters  $\Omega=5\,\rm{km}\,\rm{s}^{-1}\,
\rm{kpc}^{-1}$,$R_0=14.14$ kpc, $v_0=200\,\rm{km}\,\rm{s}^{-1}$, $\p=0.75$, and $\q=0.65$. 
{\sl Right panel:} Corresponding rotation curve. }
\label{fig:log}
\end{center}
\end{figure}

We consider a galactic model that fulfils the regularity requirements needed for the 
reduction to the centre manifold and that describes a triaxial system. We have selected the 
logarithmic potential \cite{bin87}, which is analytical and has the expression 
\begin{equation}\label{eq:log}
\Phi(x,y,z)=\frac{1}{2}v_{0}^{2}\log \left(R_{0}^{2}+x^2+\frac{y^2}{\p^2}+\frac{z^2}{\q^2}\right).
\end{equation}
The parameters $\p$ and $\q$ are non-dimensional constants that determine the shape of 
the potential, $\p$ being the planar and $\q$ the vertical axial ratios. Without loss of
generality, we can consider the potential orientated such that $0<\q<\p<1$. On the $(x,y)$ plane, 
higher values of $\p$ represent more circular systems, independent of $\q$. Analogously, on the 
$(x,z)$ plane, higher values of $\q$ represent more circular systems, independent of $\p$. For a 
standard triaxial figure, we choose $\p=0.75$ and $\q=0.65$. The constants $v_0$ and $R_0$ simply
set the scale for the velocity and length, respectively. The velocity $v_0$ corresponds to a
circular asymptotic velocity at infinity and we consider $v_0=200\,\rm{km}\,\rm{s}^{-1}$, while 
when $R_0$ is different from zero, the singularity of the potential in the origin is avoided. We choose 
$R_0=14.14\, \rm{kpc}$. The pattern speed is set to $\Omega=5\,\rm{km}\,\rm{s}^{-1}\, \rm{kpc}^{-1}$, 
characteristic for slow rotating systems \cite{bin87}. The density associated with the potential is 
obtained via Poisson's equation $\nabla^2\Phi(x,y,z)=4\pi G \rho(x,y,z)$. The terms in the Laplacian 
equation associated with the $x$ and $y$ components are always positive, while the term associated with 
the $z$-component can become negative for large values of $z$. This introduces a restriction on the 
parameters $\p$ and $\q$ of the form $\q^2 > \p^2/(1+\p^2)$. Throughout this paper we will use these 
standard values to fix the model and we will refer to it as ``Model 1''. In the left panel of 
Fig~\ref{fig:log}, we show the isodensity contours of ``Model 1'', and, in the right panel, the 
corresponding circular velocity curve, i.e. the velocity of a hypothetical star on a circular orbit. 
We note that we obtain elliptical isodensity contours and a flat rotation curve, these characteristics 
being appropriate for the study of triaxial systems.

\section{Equations of motion and study of the linear stability of the unstable
equilibrium points}
\label{sec:mot}
As is usual in galactic models, we consider a Hamiltonian formulation that, in
our case, has a Hamiltonian function given by:
$$
H(x,y,z,p_x,p_y,p_z)=\frac{1}{2}(p_x^2+p_y^2+p_z^2)+\Phi(x,y,z)-\Omega(xp_y-yp_x)\equiv E_J,
$$
where $(x,y,z)$ are the coordinate positions, $(p_x,p_y,p_z)$ are the associated conjugate 
momenta in a rotating reference frame, and $\Omega$ is the angular velocity of the galaxy 
in the inertial system (for a more detailed explanation see Binney \& Tremaine \cite{bin87}). 
The equations of motion in Hamiltonian coordinates are
$$
\begin{array}{rcl}
\dot x &=&p_x+\Omega y\\
\dot y &=&p_y-\Omega x\\
\dot z &=&p_z
\end{array} \qquad
\begin{array}{rcl}
\dot p_x&=&-\Phi_x+\Omega p_y\\
\dot p_y&=&-\Phi_y-\Omega p_x\\
\dot p_z&=&-\Phi_z,
\end{array}
$$
where $\Phi_x=\frac{\partial \Phi}{\partial x}$, $\Phi_y=\frac{\partial \Phi}{\partial y}$, and 
$\Phi_z=\frac{\partial \Phi}{\partial z}$. The effective potential is defined as $\Phi_{\hbox{\scriptsize
eff}}(x,y,z)=\Phi(x,y,z)-\Omega(xp_y-yp_x)$. The surface $\Phi_{\hbox{\scriptsize eff}}=E_J$ is called
the zero velocity surface, and its cut with the $z=0$ plane is the zero velocity curve. All regions in which 
$\Phi_{\hbox{\scriptsize eff}}>E_J$ are forbidden to a star, so we call them forbidden regions (see 
Fig.~\ref{fig:bottleneck}).

As is well-known, the effective potential associated with the galactic model has five critical points, 
located in the $xy$-plane, where 
$$\frac{\partial \Phi_{\hbox{\scriptsize eff}}}{\partial
x}=\frac{\partial \Phi_{\hbox{\scriptsize eff}}}{\partial
y}=\frac{\partial \Phi_{\hbox{\scriptsize eff}}}{\partial z}=0.$$
Due to their similarity to the corresponding points in the \rtbp, they are often called 
Lagrangian points. $L_1$ and $L_2$ are located on the $x$-axis symmetrically with respect to the
centre, $L_3$ is placed at the origin of coordinates, and $L_4$ and $L_5$ are located on the 
$y$-axis symmetrically with respect to the centre. The linearised motion around $L_1$ and $L_2$ is 
characterised by the superposition of a saddle behaviour in the $xy$-plane and 
two oscillations, one in-plane and one out-of-plane. Therefore, $L_1$ and $L_2$ are 
linearly unstable points and their behaviour is known as saddle $\times$ centre $\times$ 
centre. They are usually called the hyperbolic points. The linearised motion around $L_3$, $L_4$, 
and $L_5$ is characterised by the superposition of three oscillations, two in-plane and one 
out-of-plane. This behaviour is also known as a centre $\times$ centre $\times$ centre behaviour, 
so they are linearly stable and they are usually called the elliptic points.  

\subsection{Variation of the free parameters}
\label{sec:var}
As previously mentioned, we concentrate on energy levels slightly larger than the energies of the 
equilibrium points $L_1$ and $L_2$. At these energies there is a bottleneck in the region between 
the two zero velocity curves and we want to know the type of orbits that transit through it. 

The model parameters influence the global morphology, even for energy levels for which the bottleneck 
is open. Here we study the effect of the variation of the free parameters on the linear behaviour of 
the equilibrium points $L_1$ and $L_2$. For this purpose, we fix the values of $R_0=14.14\,\rm{kpc}$  
and $v_0=200\,\rm{km}\,\rm{s}^{-1}$, as in ``model 1'', and study the influence of $\Omega$, $\p$, 
and $\q$, separately. We make families of models in which only one of the parameters is varied, while 
the others are kept fixed. For each of the models in the family, we study the linear behaviour of the 
equilibrium points $L_1(L_2)$. We expand the effective potential, $\Phi_{\hbox{\scriptsize eff}}$, 
around one of these points and we retain only first order terms. For $L_1(L_2)$, the eigenvalues of 
the differential matrix corresponding to the planar motion are of the form: $\pm \lambda$, 
$\pm \omega\, i$, and those corresponding to the vertical motion of the form $\pm \nu\, i$. Note 
that $\lambda$, $\omega$ and $\nu$ are positive in general non-equal real numbers. Therefore, 
$L_1(L_2)$ is a linearly unstable saddle point, since the two real eigenvalues are related to the 
hyperbolic behaviour, while the purely imaginary are associated with the elliptic motion. 

\begin{figure}[!ht]
\begin{center}
\includegraphics[width=0.23\textwidth,angle=-90]{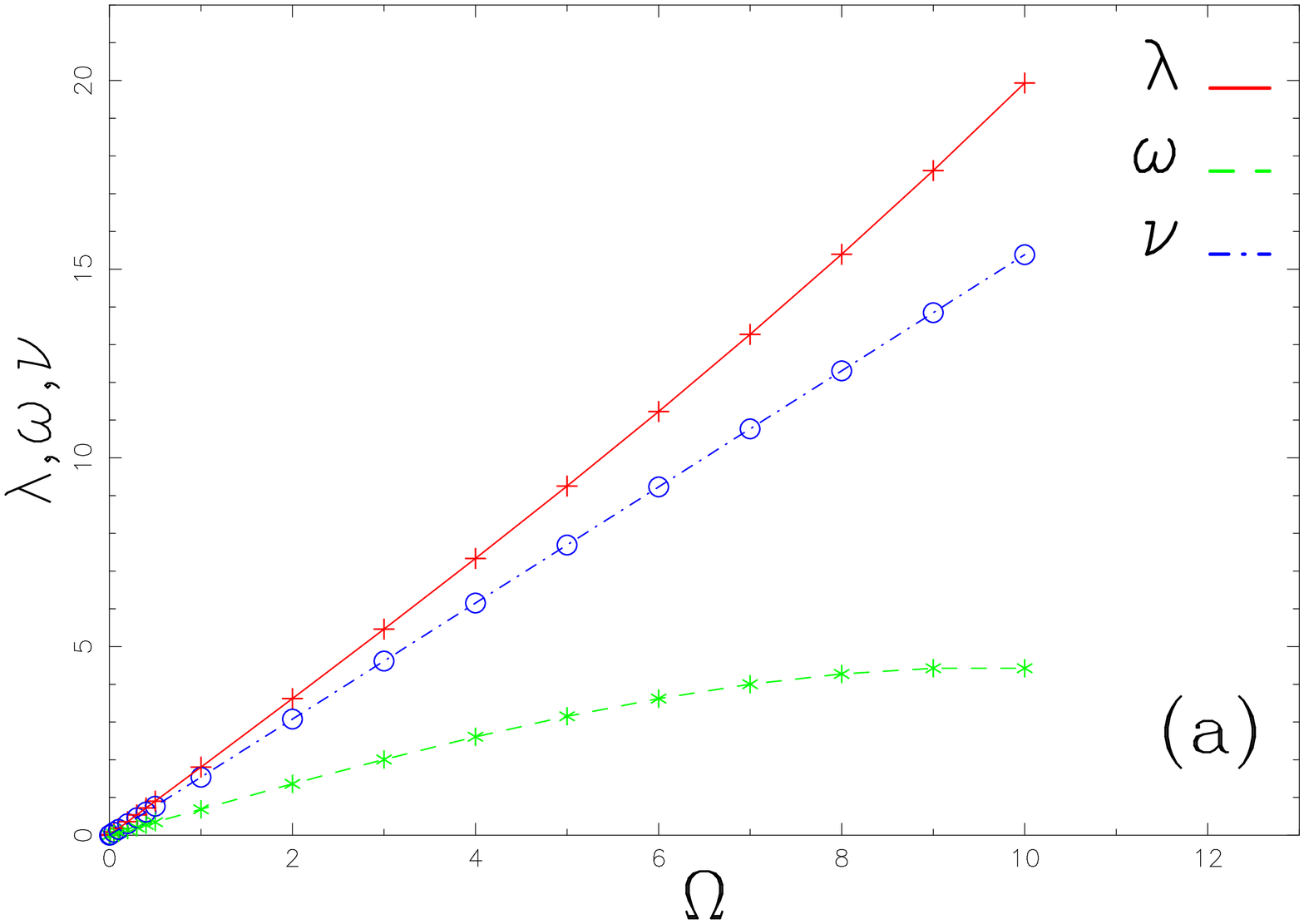}
\includegraphics[width=0.23\textwidth,angle=-90]{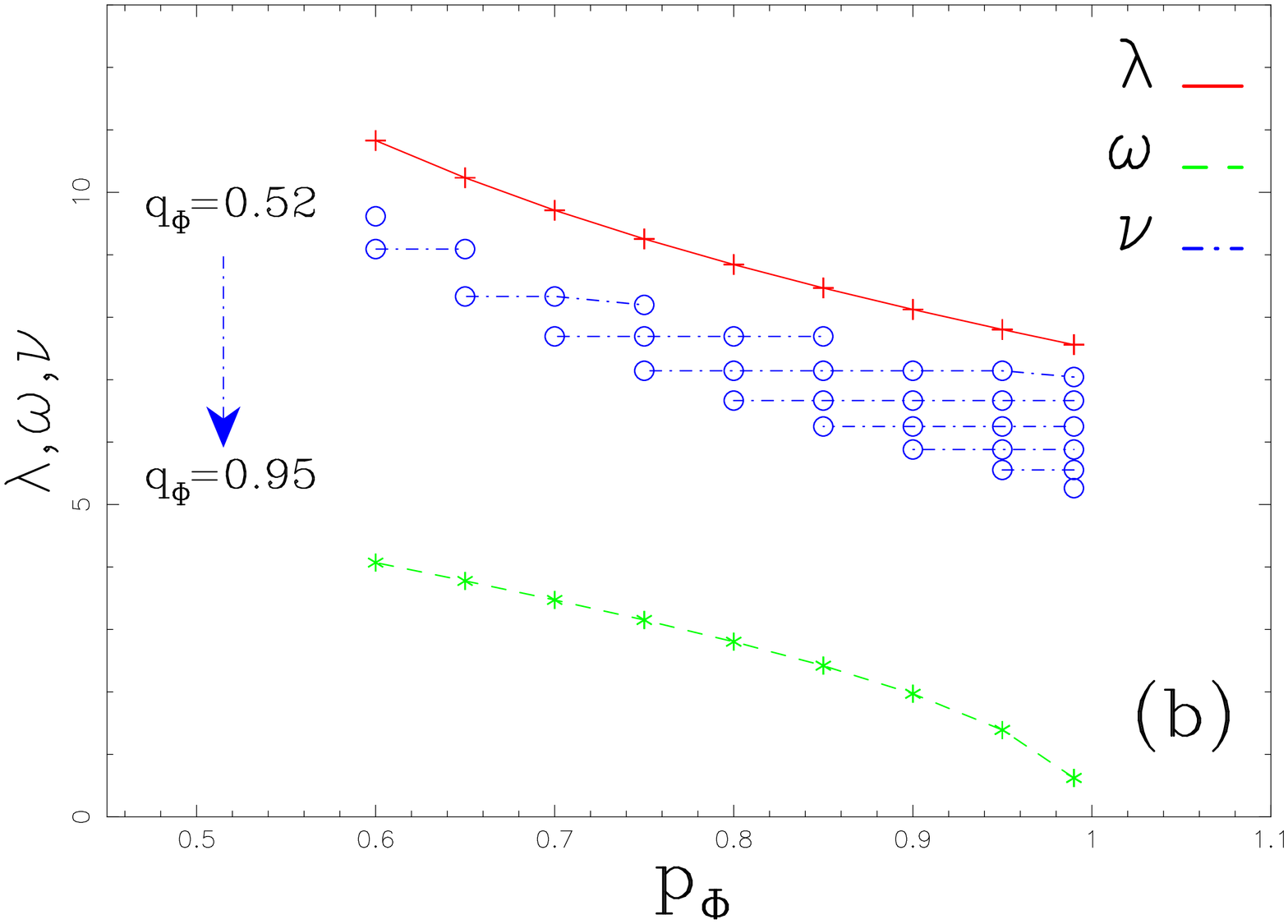}
\includegraphics[width=0.23\textwidth,angle=-90]{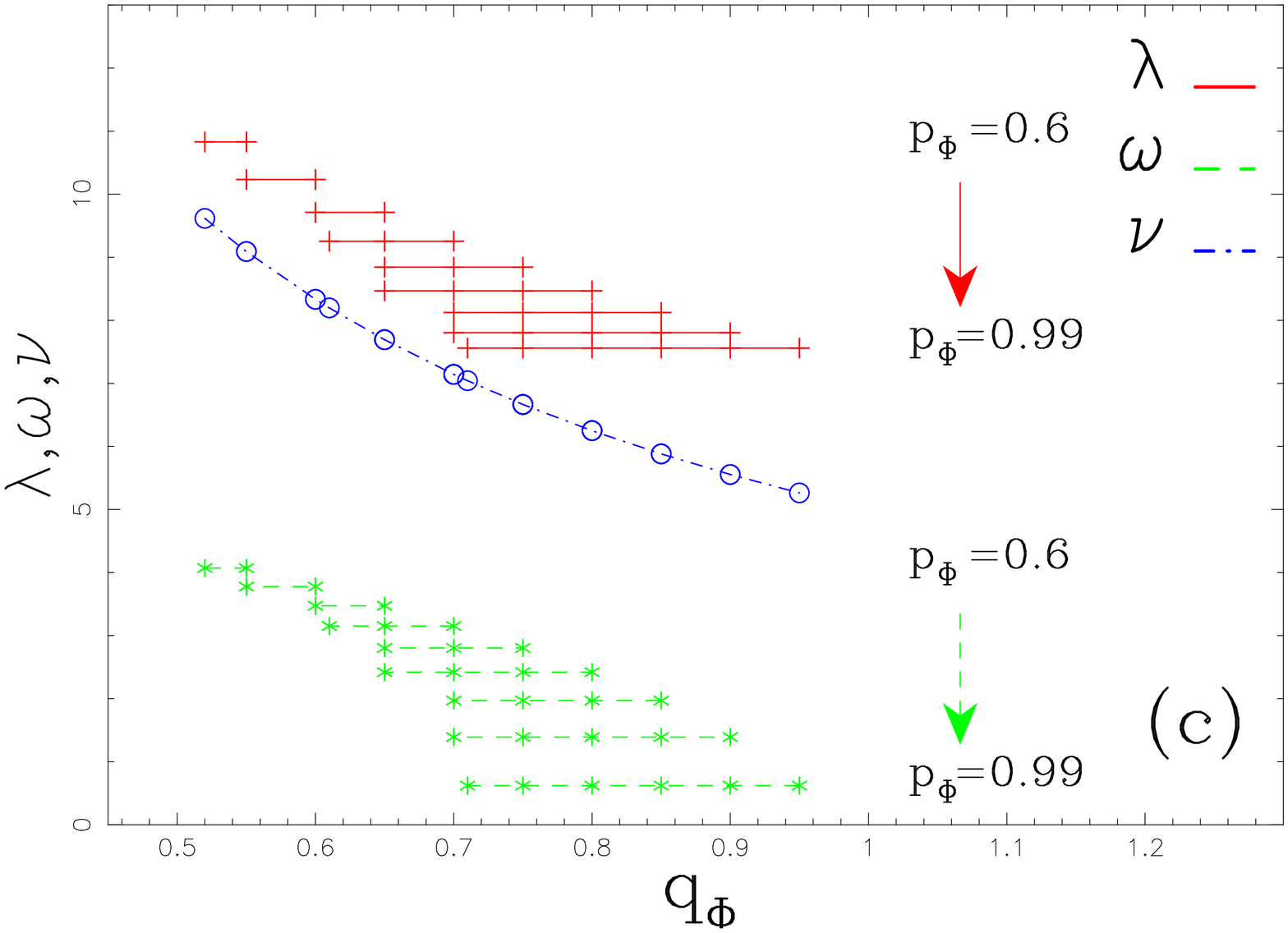}
\caption{Effect of the free parameters on the eigenvalues of the differential matrix: $\lambda$ 
(in red solid line and crosses), $\omega$ (in green dashed line and stars) and $\gamma$ (in blue dot
dashed line and open circles. The rest of the parameters are kept fixed as in ``Model 1''. {\bf (a)} 
Variation of the pattern speed, $\Omega$. {\bf (b)} Variation of the planar axial ratio, $\p$. 
{\bf (c)} Variation of the vertical axial ratio, $\q$. }
\label{fig:varvaps}
\end{center}
\end{figure}

We first study the effect of the pattern speed, $\Omega$. We vary the value of $\Omega$ within 
the range $\Omega=0.001-10\,\rm{km}\,\rm{s}^{-1}\,\rm{kpc}^{-1}$, the low values being 
characteristic for slow rotators and the large values, characteristic of fast rotators. In 
Fig.~\ref{fig:varvaps}(a), we show how the modulus of the eigenvalues changes with the 
pattern speed. We observe that as the pattern speed decreases, the eigenvalues tend to
zero, thus cancelling the hyperbolic character of $L_1(L_2)$ and removing the possibility of 
having transit orbits. The pattern speed is, on the other hand, related to the position of the 
equilibrium points through the expression $\Omega^2=r_L\left(\frac{\partial \Phi (r)}
{\partial r}\right)_{r_L}$, where $r_L$ is the distance from the centre to $L_1$ and $\Phi(r)$ is 
the potential on the equatorial plane. As $\Omega$ decreases, $r_L$ increases, so the equilibrium 
points move farther out from the centre. In Romero-G\'omez et al. \cite{rom06,rom07}, we consider a 
barred galaxy model, whose hyperbolic equilibrium points are in the vicinity of the bar ends, and we 
relate the spiral arms and rings in barred galaxies to the invariant manifolds of the periodic orbits
around the hyperbolic equilibrium points. Elliptical galaxies, however, are considered to rotate 
very slowly or even not rotate \cite{bin87,sch82}. These triaxial systems have the equilibrium points too 
farther out and, according to our results, their hyperbolic character is cancelled. Thus, no 
transfer or escape of matter is possible, which is in agreement with observations because elliptical
galaxies do not present spiral arms or rings. 

In Fig.~\ref{fig:varvaps}(b) we show the effect of the planar shape parameter, $\p$, on the 
eigenvalues of the differential matrix. We consider models where we vary $\p$ in the range 
$0.6-0.99$. We consider several models with different values of $\q$, limiting ourselves, for each 
value of $\p$, to models with a $\q$ value satisfying the restriction  $\q^2 > \p^2/(1+\p^2)$, in 
order to avoid negative densities in the $z-$axis mentioned before. The value of $\Omega$ is kept 
fixed as in ``Model 1''. In all cases, as $\p$ increases, the eigenvalues tend to decrease. In 
particular, the real eigenvalue, $\lambda$, tends to zero thus cancelling the hyperbolic character 
of the $L_1(L_2)$ points. However, note that, as expected, for a fixed value of $\q$ (vertical 
axial ratio), the corresponding vertical eigenvalue, $\nu$, does not change as we vary $\p$ 
(planar axial ratio). A similar behaviour can be observed in Fig.~\ref{fig:varvaps}(c), where we 
study the effect of the vertical shape parameter, $\q$. The eigenvalues tend to decrease as the 
vertical axial ratio increases. However, for a fixed value of the planar axial ratio, $\p$, the 
values of $\lambda$ and $\omega$ remain constant. Therefore, we can conclude that only the 
$\p$ shape parameter has an influence on the planar structure of the galaxy, while the 
parameter $\q$ will influence the vertical behaviour.

In conclusion, the hyperbolic character of $L_1$ and $L_2$ is, as expected, weaker for
slowly rotating systems or for near-axisymmetric ones, e.g elliptical galaxies. The hyperbolic 
character is strong and the bottleneck is open to global structures, in the opposite cases, e.g. bars.

\section{Reduction to the centre manifold}
\label{sec:red}
The reduction to the centre manifold is a suitable tool for studying the neighbourhood of the 
hyperbolic equilibrium points. It was first introduced in G\'omez et al. \cite{gom91} 
for the \rtbpe and further detailed in G\'omez et al. \cite{gom01}. The procedure can be summarised 
in two steps. First, the Hamiltonian function is expanded in power series around the 
equilibrium point. Then, a partial normal form scheme is applied in order to uncouple (up to a high 
order) the hyperbolic directions from the elliptic ones. Now the truncated Hamiltonian has an 
invariant manifold tangent to the elliptic directions of the linear part. The restriction to the 
invariant manifold tangent to the elliptic directions leads to a Hamiltonian system with two degrees 
of freedom and an elliptic equilibrium point at the origin. This restriction to the manifold is the 
so-called reduction to the centre manifold and the study of the dynamics of the reduced Hamiltonian 
gives a qualitative description of the phase space near the equilibrium point. 

We want to note here that we do not use the Birkhoff Normal Form (hereafter, BNF). We are
interested in performing a local study of a hyperbolic equilibrium point in the largest region
possible. In that sense, the reduction to the centre manifold can provide better results in a 
larger neighbourhood than BNF. Furthermore, as we will see in Section~\ref{sec:res}, we are 
interested in computing the invariant manifolds associated with periodic orbits and quasi-periodic 
orbits and the normal form described here provides as well initial conditions to compute both invariant 
objects.

We proceed by first obtaining the quadratic real normal form of the second order term of 
the Hamiltonian function. This is achieved by performing two changes of coordinates. The first 
one consists of a translation plus an homothecy in order to place the new origin of coordinates 
on the equilibrium point and to set the unit length equal to the distance from the equilibrium
point to the origin. The second change of variables allows us to write the second order part of the 
Hamiltonian in its real normal form. We then perform a complexification of the variables, in 
order to obtain a Hamiltonian of diagonal form. Finally, the normal form of higher order terms
is computed. In the following subsections, we describe in detail the steps performed and in 
the last subsection we study the practical convergence of the reduced Hamiltonian.

\subsection{Translation and homothecy}

As previously mentioned, the first change of coordinates consists of a translation plus 
an homothecy. The translation sets the new origin of coordinates on any of the symmetric 
hyperbolic equilibrium points and the homothecy sets the unit of distance as the distance from the 
origin of the old system to the position of the equilibrium point. Furthermore, the change has 
to be symplectic, i.e. it must preserve the Hamiltonian equations in the new variables.

If we write $(x^\prime,y^\prime,z^\prime,p_x^\prime,p_y^\prime,p_z^\prime)$ as the
Hamiltonian barycentric coordinates and \\$(x,y,z,p_x,p_y,p_z)$ as the
Hamiltonian ones around the hyperbolic equilibrium point, the symplectic change can be
written as:
$$
\begin{array}{rcl}
x &=& ( x^\prime -x_{L_i})/\gamma\rule[-.5cm]{0cm}{0.5cm}\\
p_x&=&\gamma(p_x^\prime-p_{x_{L_i}})\rule[-.5cm]{0cm}{0.5cm}
\end{array} \qquad
\begin{array}{rcl}
y &=& ( y^\prime -y_{L_i})/\gamma\rule[-.5cm]{0cm}{0.5cm}\\
p_y&=&\gamma(p_y^\prime-p_{y_{L_i}})\rule[-.5cm]{0cm}{0.5cm}\\
\end{array}
\qquad \begin{array}{rcl}
z &=& z^\prime/\gamma\rule[-.5cm]{0cm}{0.5cm}\\
p_z&=&\gamma p_z^\prime\rule[-.5cm]{0cm}{0.5cm},
\end{array}
$$
where $L_i=(x_{L_i},y_{L_i},0,p_{x_{L_i}},p_{y_{L_i}},0),\, i=1,2$ are the coordinates of the hyperbolic 
equilibrium point $L_i$ and $\gamma = \sqrt{x_{L_i}^2+y_{L_i}^2}$.
Expanding the logarithmic potential up to second order, the expression for the Hamiltonian
in the new coordinates around $L_i,\, i=1,2$, is
\begin{eqnarray*}
H(x,y,z,p_x,p_y,p_z)&=&H_0+\frac{1}{2}\frac{1}{\gamma^2}(p_x^2+p_y^2+p_z^2)+\frac{1}{2}v_0^2 \gamma^2\left( \left(\frac{1}{K}-
\frac{2x_{L_i}^2}{K^2}\right)x^2+\right.\\
   & & \left. +\left(\frac{1}{K\p^2}-\frac{2y_{L_i}^2}{K^2\p^4}\right)y^2-\frac{4x_{L_i}y_{L_i}}{K^2\p^2}xy+\frac{1}{K\q^2}z^2\right)-\\
   & & -\Omega (xp_y+yp_x)+\ldots
\end{eqnarray*}
where $H_0=\frac{1}{2}v_0^2\ln K +\frac{1}{2}(p_{x_{L_i}}^2+p_{y_{L_i}}^2)+\Omega (y_{L_i}p_{x_{L_i}}- x_{L_i}p_{y_{L_i}})$
is a constant and $K=R_0^2+x_{L_i}^2+\frac{y_{L_i}^2}{\p^2}$.  
 
\subsection{The quadratic real normal form}

The second change of coordinates will express the second order term of
the Hamiltonian in its real normal form 
$H_2=\lambda x p_x+\frac{1}{2}\omega_1(y^2+p_y^2)+\frac{1}{2}\omega_2(z^2+p_z^2)$.
This is accomplished via a linear and symplectic change performed
around the hyperbolic equilibrium point $L_1$. The analysis for $L_2$ is analogous and the results
are symmetric.

For $L_1$, we have
$$x_{L_1}=\sqrt{\left( \frac{v_0^2}{\Omega^2}\right)-R_0^2},\quad y_{L_1}=0,\quad p_{x_{L_1}}=0,\quad p_{y_{L_1}}=\Omega \sqrt{\left( \frac{v_0^2}{\Omega^2}\right)-R_0^2}$$
and $\gamma=x_{L_1}$. 

We note that the planar motion is uncoupled from the vertical motion, so we perform the reduction
to the real normal form in two steps, that is, first regarding the in-plane motion and then the
out-of-plane motion. In the 2D case, the differential matrix around $L_1$ becomes
$$
M=\left (
\begin{array}{cccc}
    0 & a & b & 0 \\
   -a & 0 & 0 & b \\
    c & 0 & 0 & a \\
    0 & -d & -a & 0
\end{array} \right ),
$$
where
$$a=\Omega,\quad b=\frac{1}{\gamma^2},\quad c=\frac{\gamma^2}{v_0^2}\left( \frac{v_0^2}{\Omega^2}-2R_0^2\right)\Omega^4,\quad d=\frac{\Omega^2 \gamma^2}{\p^2}.$$
The matrix $M$ has two real, $\pm \lambda$, and two purely imaginary,
$\pm \omega_1 i$, eigenvalues. The eigenvectors associated with the real eigenvalues have the 
following expression
\begin{equation}\label{eq:vec}
v_{\lambda}= \left[ \begin{array}{c}
                    2\frac{\lambda}{b}a \\
		    -\left( c-\frac{\lambda^2}{b}+\frac{a^2}{b}\right)\\
		    2\frac{\lambda^2}{b^2}a -\frac{a}{b}\left(-c+\frac{\lambda^2}{b}-\frac{a^2}{b}\right)\\
		    \frac{2}{b^2}a^2\lambda+\frac{\lambda}{b}\left(-c+\frac{\lambda^2}{b}-\frac{a^2}{b}\right)
		  \end{array} \right],
\end{equation}
while the eigenvector associated with the imaginary eigenvalue is 
$$
v_\lambda= u+iv=\left[ \begin{array}{c}
                    0 \\
		    -c-\frac{\omega_1^2}{b}-\frac{a^2}{b}\\
		    -2\frac{\omega_1^2}{b^2}a -\frac{a}{b}\left(-c+\frac{\omega_1^2}{b}-\frac{a^2}{b}\right)\\
		    0
		  \end{array} \right] +i
                \left[ \begin{array}{c}
		    2\frac{a}{b}\omega_1 \\
		    0\\
		    0\\
		    2\frac{a^2}{b^2}\omega_1+\frac{\omega_1}{b}\left(-c-\frac{\omega_1^2}{b}-\frac{a^2}{b}\right)
		  \end{array} \right]
$$
Considering $V^{{\prime}^T}\,J\,V^{\prime}$, where  
$V^{\prime}=\left(v_{+\lambda}\, u \,\, v_{-\lambda} \, v \right)$, we obtain the suitable
scaling factors $d_{\lambda}$ and $d_{\omega_1}$ that give the final symplectic change of variables 
regarding the in-plane motion:
$$V=\left(\frac{v_{+\lambda}}{\sqrt{d_{\lambda}}} \quad \frac{u}{\sqrt{d_{\omega_1}}} \quad 
\frac{v_{-\lambda}}{\sqrt{d_{\lambda}}} \quad \frac{v}{\sqrt{d_{\omega_1}}} \right).$$ 

Regarding the out-of-plane motion, the symplectic transformation is given by
$$
z \rightarrow \frac{\gamma}{\sqrt{\omega_2}}z \qquad p_z \rightarrow \frac{\sqrt{\omega_2}}{\gamma}p_z,
$$
where $\omega_2^2=\frac{v_0^2}{K\q^2}>0$ always. The final change of coordinates
sets the Hamiltonian around the equilibrium point $L_1$
with coordinates that express the second order term in its real normal
form. That is, after performing the second change of
coordinates, the second order term of the Hamiltonian is:
\begin{equation}\label{eq:ham2}
H_2=\lambda
xp_x+\frac{\omega_1}{2}\left(y^2+p_y^2\right)+\frac{\omega_2}{2}\left(z^2+p_z^2\right).
\end{equation}

For subsequent computational purposes, it is desirable to work with
complex variables. Then, a third change of variables is needed in
order to write the second order term of the Hamiltonian in the complex
diagonal normal form. 

\subsection{Complexification}
We are interested in writing the second order Hamiltonian in a diagonal form, in order to 
simplify the homologic equation that we encounter in the normal form process. We perform a 
complexification of the variables related to $y$ and $z$. That is, we perform a symplectic change
of coordinates from real coordinates $(x,y,z,p_x,p_y,p_z)$ to complex coordinates of the normal 
form $(q_1,q_2,q_3,p_1,p_2,p_3)$ complexifying the pairs $(y,p_y)$ and $(z,p_z)$ by means of:
$$
\begin{array}{rcl}
q_1 &=&x\rule[-.5cm]{0cm}{0.5cm}\\
p_1&=& p_x \rule[-.5cm]{0cm}{0.5cm}
\end{array} \qquad
\begin{array}{rcl}
q_2 &=&\displaystyle\frac{y-ip_y}{\sqrt{2}}\rule[-.5cm]{0cm}{0.5cm}\\
p_2&=&\displaystyle\frac{p_y-iy}{\sqrt{2}}\rule[-.5cm]{0cm}{0.5cm}\\
\end{array}\qquad
\begin{array}{rcl}
q_3 &=& \displaystyle\frac{z-ip_z}{\sqrt{2}}\rule[-.5cm]{0cm}{0.5cm}\\
p_3&=&\displaystyle\frac{p_z-iz}{\sqrt{2}}\rule[-.5cm]{0cm}{0.5cm}.
\end{array}
$$
In this way, the second order Hamiltonian has the form:
\begin{equation}\label{eq:cpx}
H_2=\lambda q_1p_1+i\omega_1q_2p_2+i\omega_2q_3p_3,
\end{equation}
where $\lambda$, $\omega_1$ and $\omega_2$ are real positive numbers.

\subsection{Normal form of higher order terms}
We have now written the second order term of the Hamiltonian, $H_2$,
in its real and in its complex normal form. The final step to obtain the
centre manifold consists of removing some monomials in the expansion
of the Hamiltonian so that the final Hamiltonian has an invariant
manifold tangent to the elliptic directions of $H_2$. We apply the method 
known as the Lie Series method (see \cite{jor99,gom01} and references therein for a 
detailed description in the \rtbp). Until the end of this section, we use the following notation. 
If $x=(x_1,\dots,x_n)$  is a vector of complex numbers and $k=(k_1,\dots,k_n)$ is an integer vector, 
we denote $x^k$ the term $x_1^{k_1}\dots x_n^{k_n}$ and in this context we define $0^0$ as $1$. 
We define $\vert k\vert$ as $\sum_j\vert k_j \vert$. 

The initial Hamiltonian is expanded around the hyperbolic equilibrium
point in the complex coordinates for which the second order term is in
diagonal form (Eq. \ref{eq:cpx}). Thus, the Hamiltonian\footnote{Bold letters 
denote vector notation.} has the form 
$$
H({\bf q},{\bf p})=H_2({\bf q},{\bf p})+\displaystyle \sum_{n\ge 3} H_n({\bf q},{\bf p}),
$$
where $H_2$ is given as in Eq. \ref{eq:cpx} and $H_n$ is a
homogeneous polynomial of degree $n$ of the form $\displaystyle
\sum_{i,j} h_{ij}q_1^{i_1}p_1^{j_1}q_2^{i_2}p_2^{j_2}q_3^{i_3}p_3^{j_3}$, where
$h_{ij}$ denotes $h_{i_1,i_2,i_3,j_1,j_2,j_3}$.

The Poisson bracket of two functions, $F({\bf q},{\bf p})$ and $G({\bf q},{\bf p})$, of position 
and momenta is defined as
$$
\{F,G\}=\displaystyle \sum_{i=1}^3\left (\frac{\partial F}{\partial q_i}\frac{\partial G}{\partial p_i}-
\frac{\partial F}{\partial p_i}\frac{\partial G}{\partial q_i} \right).
$$
If $F$ and $G$ are homogeneous polynomials of degree $r$ and $s$ respectively,
then $\{F,G\}$ is a homogeneous polynomial of degree $r+s-2$. If a function $G({\bf q},{\bf p})$ is
considered a Hamiltonian then, the function $\hat H$ defined by
$$
\hat H \equiv H+\{H,G\}+\frac{1}{2!}\{\{H,G\},G\}+\frac{1}{3!}\{\{\{H,G\},G\},G\}+...,
$$
is the result of applying a time one map canonical transformation under the flow of $G$ to the 
Hamiltonian $H$. The Hamiltonian $G$ is usually called the generating function. If we choose a 
homogeneous polynomial of degree 3, $G_3$, as a generating function, then the homogeneous 
polynomials of degree $n$, $\hat H_n$, such that $\hat H=\sum_{n\ge 2}\hat H_n$ are given by

$$
\begin{array}{rcl}
\hat H_2 & = & H_2,\\
\hat H_3 & = & H_3+\{H_2,G_3\},\\
\hat H_4 & = & H_4+\{H_3,G_3\}+\frac{1}{2!}\{\{H_2,G_3\},G_3\},\\
\vdots & &
\end{array}
$$
In order to remove all the terms of order three in the new Hamiltonian, i.e. to have $\hat H_3=0$, 
we must choose a generating function, $G_3$, such that it solves the homological equation 
$\{H_2,G_3\}=-H_3$. This procedure can be used recurrently trying to find homogeneous polynomials 
$G_n$ to remove non-resonant terms of the Hamiltonian.

In our case, we are interested in removing the instability associated with the hyperbolic character 
of the Hamiltonian $H$. We note that the second order term of the Hamiltonian provides the linear part
of the Hamiltonian equations. Therefore, the instability is associated with the term $\lambda q_1p_1$. 
For the linear approximation of the Hamiltonian equations, the centre part can be obtained by setting
$q_1=p_1=0$. If we want the trajectory to remain tangent to this space (i.e. to have $q_1(t)=p_1(t)=0$ 
for all $t>0$), then we need $\dot q_1(0)=\dot p_1(0)=0$ when adding the nonlinear terms. Due to the
autonomous character of the Hamiltonian system, we will obtain $q_1(t)=p_1(t)=0$ for all $t\ge0$. We 
remember that the Hamiltonian equations associated with a Hamiltonian $H({\bf q},{\bf p})$ are
$$\dot q_i=\frac{\partial H}{\partial p_i};\quad \dot
p_i=-\frac{\partial H}{\partial q_i}. $$ 
In particular,
$$
\begin{array}{rcccl}
\dot q_1 & = & \frac{\partial H}{\partial p_1} & = & \lambda
q_1+\displaystyle \sum_{n\ne 3}
h_{ij}q_1^{i_1}p_1^{j_1-1}q_2^{i_2}p_2^{j_2}q_3^{i_3}p_3^{j_3}\\
\dot p_1 & = & -\frac{\partial H}{\partial q_1} & = & -\lambda
p_1-\displaystyle \sum_{n\ne 3}
h_{ij}q_1^{i_1-1}p_1^{j_1}q_2^{i_2}p_2^{j_2}q_3^{i_3}p_3^{j_3}.\\
\end{array}
$$
Therefore, we can obtain the required condition, $\dot q_1(0)=\dot p_1(0)=0$ when $q_1(0)=p_1(0)=0$, if 
in the expansion of the Hamiltonian $H$ there are no monomials with exponents $(1,0,i_2,j_2,i_3,j_3)$ 
and $(0,1,i_2,j_2,i_3,j_3)$. Different changes of variables can be used for such purpose, the most
common being the ones cancelling the terms with $i_1\ne j_1$ \cite{gom01} or cancelling only terms with
$i_1+j_1=1$. This second choice cancels the minimum number of terms in the Hamiltonian and is the
one we have chosen in this paper.

Thus, after the final change of variables the Hamiltonian can be written in the form:
$$
\overline H({\bf q},{\bf p})=\overline H_N({\bf q},{\bf p})+R_N({\bf q},{\bf p}),
$$
where $\overline H_N({\bf q},{\bf p})$ is a polynomial of degree $N$ without terms
of $i_1+j_1=1$ and $R_N$ is a remainder of order $N+1$, which is very
small near $L_1(L_2)$. 

Finally, using the inverse change of variables of the complexification, the truncated Hamiltonian, 
$\overline H_N$ can be expanded in real form and we obtain
$$
\overline H_N({\bf q},{\bf p})=H_2({\bf q},{\bf p})+\displaystyle \sum_{n=3}^{N}H_n({\bf q},{\bf p}),
$$
where the second order term $H_2({\bf q},{\bf p})$ is as in Eq. \ref{eq:ham2}.

\subsection{Practical convergence}

\begin{figure}[!ht]
\begin{center}
\includegraphics[scale=0.35,angle=-90]{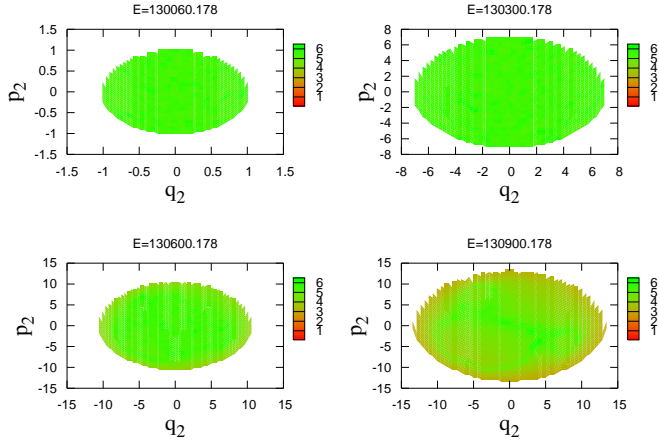}
\caption{Study on the practical convergence of the reduction to the centre manifold for a tolerance of 
$\epsilon=10^{-6}$ in ``model 1''. In all panels, we associate to each initial condition and the final 
integration time to a colour in the palette. The energy levels used are $130060.178$ (top left), 
$130300.178$ (top right), $130600.178$ (bottom left) and $130900.178$ (bottom right).}
\label{fig:convA}
\end{center}
\end{figure}

\begin{figure}[!ht]
\begin{center}
\includegraphics[scale=0.35,angle=-90]{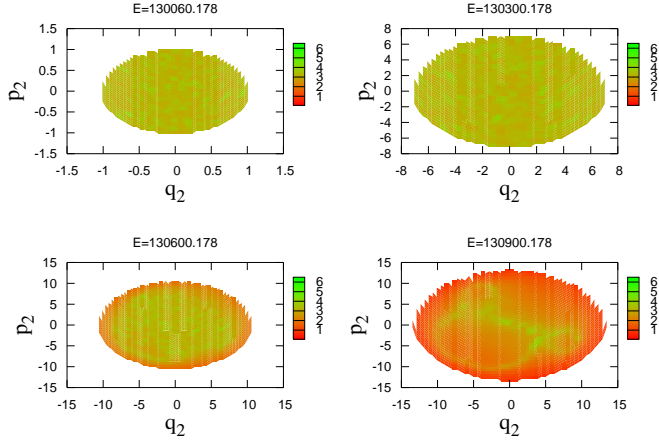}
\caption{As in Fig.~\ref{fig:convA} for a tolerance of $\epsilon=10^{-9}$}
\label{fig:convB}
\end{center}
\end{figure}

In order to check the practical convergence of the truncated series, we integrate the same initial 
condition on the centre manifold using both the reduced Hamiltonian up to order $15$ and the Hamiltonian 
in barycentric coordinates. At each time step, we compare the two position vectors by computing the norm 
of the difference. When this becomes greater than a given tolerance, we plot the initial condition and 
the time step it has reached. In Fig.~\ref{fig:convA}, we show the practical convergence plots 
for four different energy levels and a tolerance of $\epsilon=10^{-6}$ for ``model 1''. In 
Fig.~\ref{fig:convB}, we plot the energy levels for the same model but with a tolerance of $\epsilon=10^{-9}$. 
The energy of the equilibrium point is $E_J(L_1)=130055.178$. In each panel, we consider an energy 
level higher than the previous one and initial conditions on the centre manifold. We note that, due to the 
energy confinement, the initial conditions valid for a given energy level lie inside the planar Lyapunov 
orbit. As the energy increases, the time, $t$, at which the differences become greater than the tolerance, 
decreases. This is reflected in Figs.~\ref{fig:convA} and \ref{fig:convB} by the colour palette. In 
all panels, the colour palette ranges from $t=0$ to $t=6.5$, the green colour being associated with large 
times and the red colour, to small times. In Fig.~\ref{fig:convA}, where we use $\epsilon=10^{-6}$ as 
tolerance, we note that in the top left panel, the dominating colour is green, while as we increase 
the energy, the colours tend to yellow and red. We also observe that, as expected, as the tolerance 
decreases, the time also decreases, so in the bottom right panel of Fig.~\ref{fig:convB} the region is 
practically all orange and red. 

Note that for a tolerance of $10^{-6}$, the practical convergence of the normal form is good for
a wide range of energies, as all panels in Fig.~\ref{fig:convA} are basically green. Even for a more 
restrictive tolerance, $10^{-9}$, the normal form gives good results for a considerably range of energies,
see top panels of Fig.~\ref{fig:convB}. Also note that the numerical test we use here integrates initial 
conditions using both the reduced and the initial Hamiltonian. Since the initial Hamiltonian contains the 
instability, that grows exponentially in time, the time estimations we show here underestimate the real 
values of $t$ for which the reduced Hamiltonian gives accurate results. 

\section{Invariant objects around $L_1/L_2$. Transfer of matter}
\label{sec:res}
In the following, we apply the methodology introduced in the previous section to compute the invariant 
manifolds associated with both periodic orbits and invariant tori. We fix the values of the potential and 
study the neighbourhood of the hyperbolic equilibrium point $L_1$. We first compute the invariant curves 
around the equilibrium point on the section given by ${\cal I}=\{z=0\}$ and within a range of energies. 
We then compute the invariant manifolds of the planar and vertical Lyapunov orbits and of the invariant 
tori. We also study the role invariant manifolds of periodic orbits have in the transfer of
matter in the galaxy by computing the homoclinic and heteroclinic connections between the
periodic orbits of a given energy level. Finally, we study the variation of the free 
parameters and check their influence on the presence of homoclinic and heteroclinic orbits,
thus, on the global shape produced by the invariant manifolds.

\subsection{Invariant manifolds of periodic orbits and quasi-periodic
orbits}
\label{sec:invman}
Throughout this section we fix the model parameters to those of the standard ``model 1'' and
we compute the invariant objects using the Hamiltonian in the normal form up to order $15$ in the 
reduced coordinates. For clarity of the representation, however, at each step we perform the backwards 
change of variables to plot the results in the initial frame using barycentric coordinates. 

\begin{figure}[!ht]
\begin{center}
\includegraphics[scale=0.3,angle=-90]{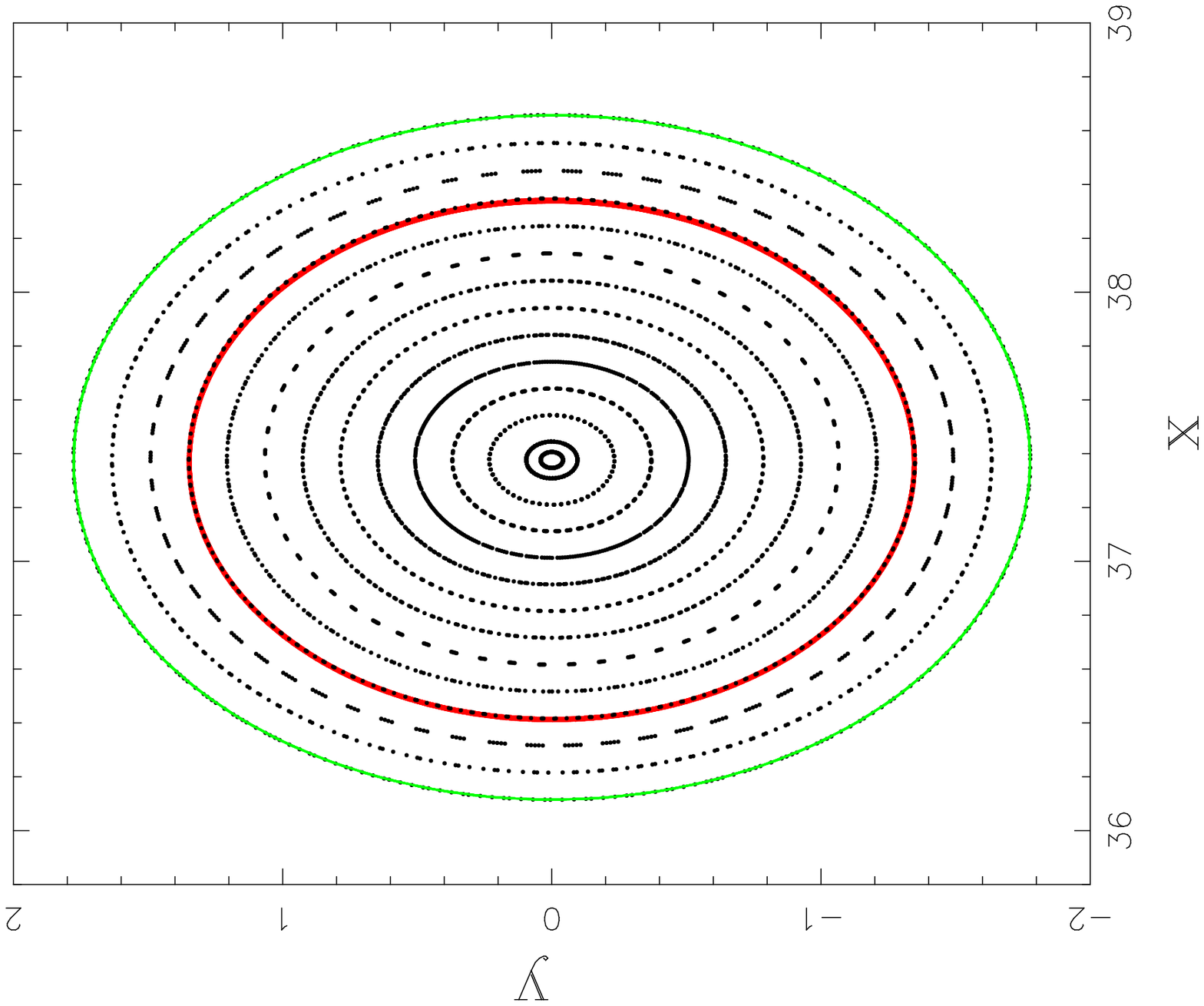} \\
\includegraphics[scale=0.22,angle=-90]{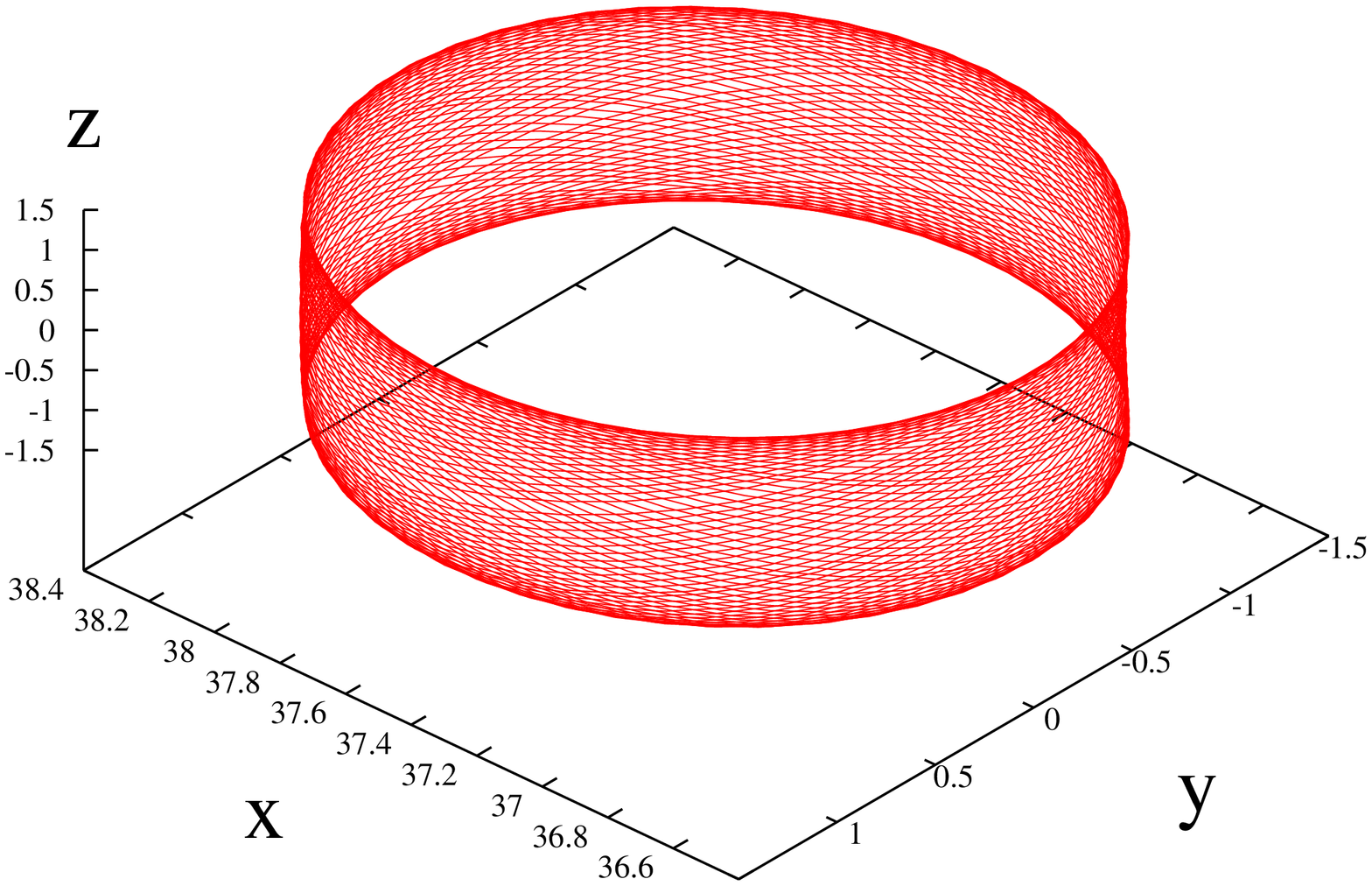}
\includegraphics[scale=0.22,angle=-90]{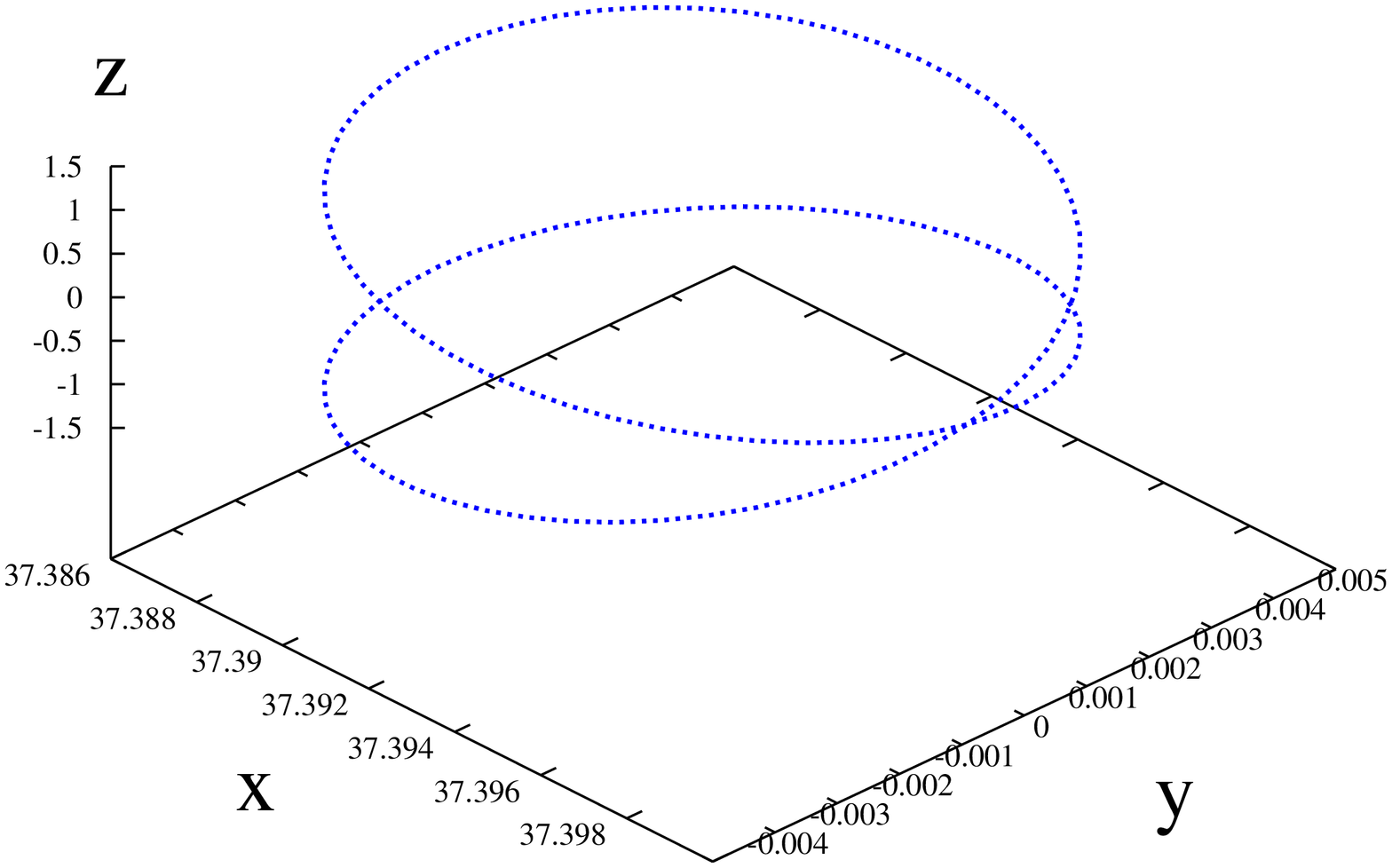}
\caption{{\sl Top panel:} Invariant curves (black dots) and planar Lyapunov orbit 
(green solid line) around the equilibrium point $L_1$ for the logarithmic potential 
with parameters as in ``model 1'', and energy level $E_J=130100.178$. The red solid line 
corresponds to the invariant torus plotted in the bottom left panel. {\sl Bottom left 
panel:} Invariant torus (red solid line) corresponding to the (red solid line) invariant 
curve on the section $\cal{I}$. {\sl Bottom right panel:} Vertical Lyapunov orbit 
(blue dotted line) of the same energy.}
\label{fig:invcor}
\end{center}
\end{figure}

We compute the invariant curves on the section $\cal{I}$ and the planar Lyapunov orbit surrounding 
the equilibrium point, see top panel of Fig.~\ref{fig:invcor}. To obtain the invariant curves, 
we take initial conditions on the galactic plane of the initial frame of reference and we plot 
the cuts with the section. Each invariant curve on the section represents an invariant torus, as is
shown in the top and bottom left panels of Fig.~\ref{fig:invcor}. In the bottom left panel, we plot 
the invariant torus corresponding to the invariant curve on the section $\cal{I}$ (red solid line) in the 
original coordinates, although the numerical integration is performed using the reduced Hamiltonian. 
Since the reduced Hamiltonian has the hyperbolic component decoupled from the elliptic one, we can 
perform long time integrations without the instability. For convenience, however, we use the normal
form change of variables at each integration step to display the orbits in the original coordinates.
In the bottom right panel, we plot the vertical Lyapunov orbit. We observe that for the range of 
energies for which the normal form of the Hamiltonian is accurate, the section $\cal{I}$ of the 
invariant curves has always the same aspect, i.e., the only macroscopic objects we see are the 
Lyapunov orbits and invariant tori.

\begin{figure}[!ht]
\begin{center}
\includegraphics[scale=0.4,angle=-90]{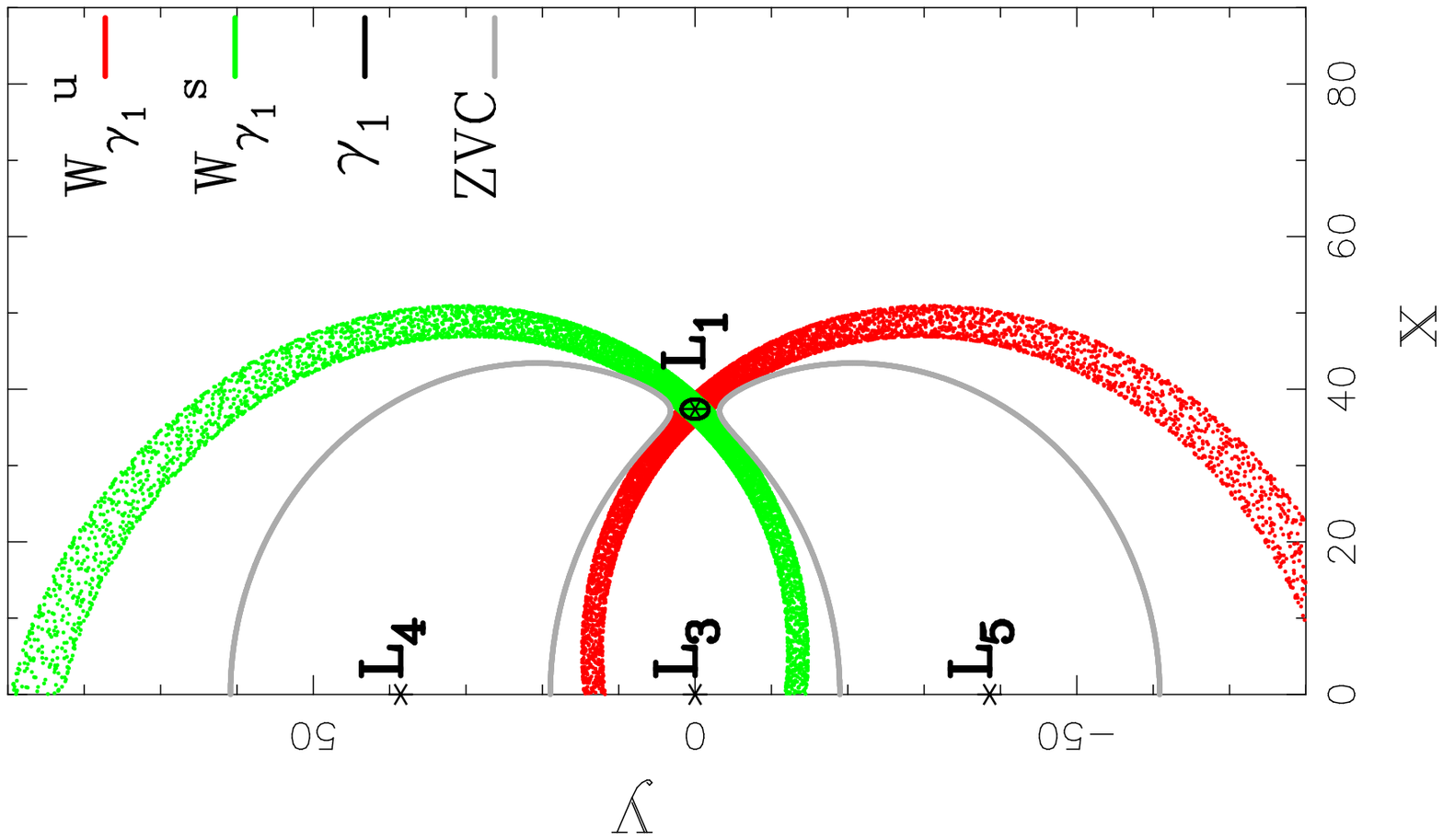}\hspace{1.5cm}
\includegraphics[scale=0.4,angle=-90]{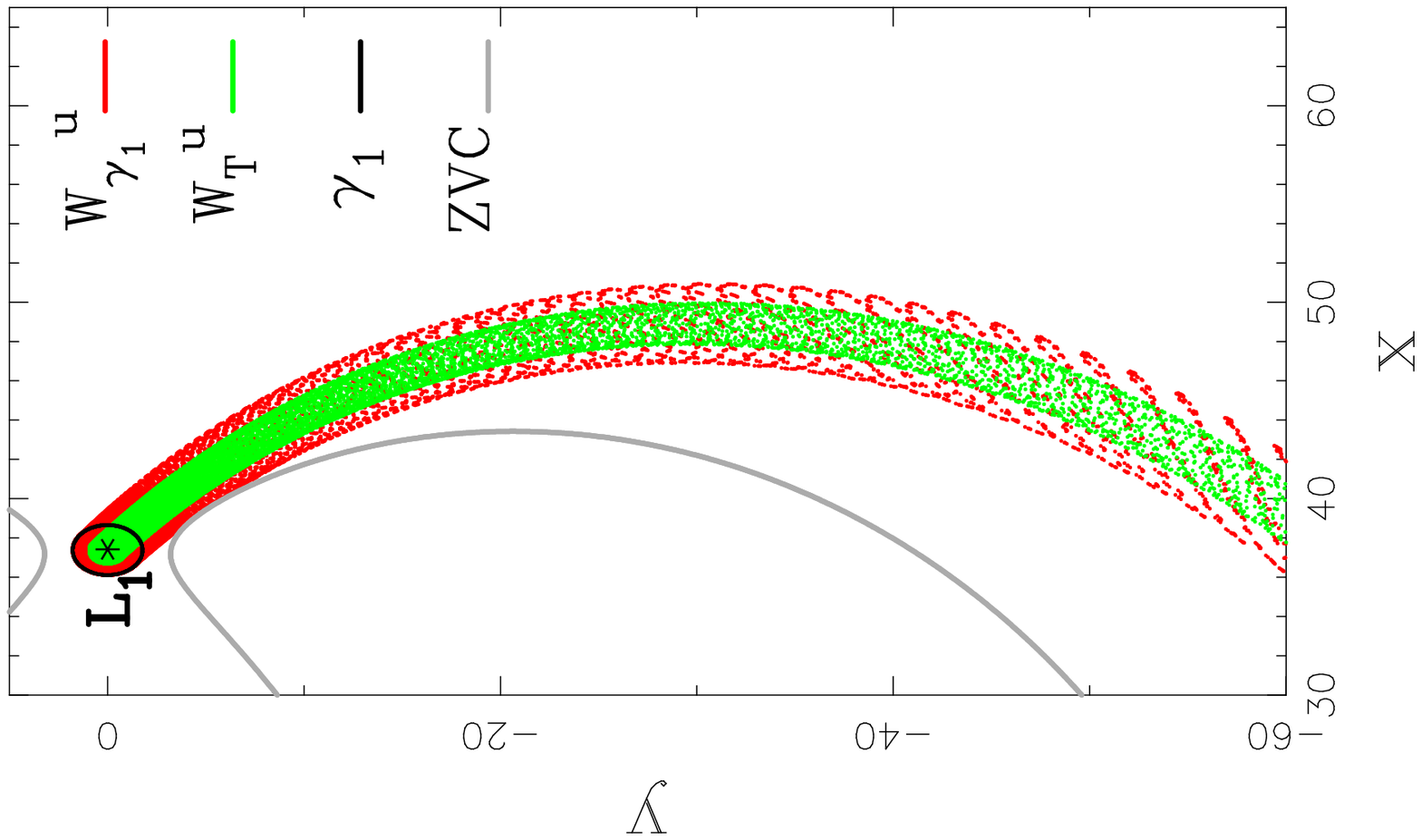}
\end{center}
\caption{Invariant manifolds of a periodic orbit (left panel) and of a quasi-periodic orbit 
(right panel) of the logarithmic model for the energy level of $E_J=130155.178$. 
{\sl Left panel:} In the centre of the plot, a black solid line shows the planar Lyapunov 
orbit around $L_1$. We plot the two branches of the unstable invariant manifold (red dotted lines), 
the two branches of the stable invariant manifold (green dotted lines), and the the zero velocity 
curves (grey solid lines) defining the forbidden region. {\sl Right panel:} Invariant manifolds 
of an invariant torus, $W^u_T$ (in green), and of a periodic orbit, $W^u_{\gamma_1}$ (in red). In 
black solid line, the planar Lyapunov orbit around $L_1$ and, in grey, the zero velocity curves.}
\label{fig:brlog}
\end{figure}

We can also compute the hyperbolic invariant manifolds associated with the planar and vertical 
Lyapunov orbits and with the invariant tori. We first compute the stable and unstable invariant manifolds
associated with the planar Lyapunov periodic orbit around the equilibrium point $L_1$, 
$W_{\gamma_1}^s$ and $W_{\gamma_1}^u$ respectively, using the truncated Hamiltonian 
up to order $15$. According to the normal form scheme used, the initial
conditions approximating the unstable invariant manifold of a planar Lyapunov orbit 
have the form $(q_1,p_1,q_2,p_2,q_3,p_3)=(\pm \epsilon,0,q_2,p_2,0,0)$ and the
initial conditions for the stable invariant manifold are $(q_1,p_1,q_2,p_2,q_3,p_3)$ 
$=(0,\pm \epsilon,q_2,p_2,0,0)$. Note that invariant manifolds are sets of asymptotic 
orbits that tend to and depart from the periodic orbit. In the left panel of Fig.~\ref{fig:brlog}, 
we plot $W_{\gamma_1}^s$ and $W_{\gamma_1}^u$ for the energy level $E_J=130155.178$ 
($E_J(L_1)=130055.178$) and $\epsilon=10^{-5}$. We also plot the zero velocity curve defining the
forbidden region. Here we observe the characteristic saddle behaviour of the invariant manifolds 
due to the hyperbolic behaviour of the periodic orbit. Also note that the invariant manifolds 
connect the inner region with the outer region, in the sense that the existence of invariant manifolds 
of periodic orbits implies also the existence of transit orbits between the two regions delimited by 
the zero velocity curves. In the right panel of Fig.~\ref{fig:brlog}, we plot 
the unstable invariant manifold of an invariant torus for the energy level $E_J=130155.178$. 
The initial conditions to compute the unstable and stable invariant manifolds of invariant
tori are easily provided by the normal form we use. That is, 
$(q_1,p_1,q_2,p_2,q_3,p_3)=(\pm \epsilon,0,q_2,p_2,q_3,p_3)$ for the unstable manifold,
$(0,\pm \epsilon,q_2,p_2,q_3,p_3)$ for the stable manifold and 
$\epsilon=10^{-5}$. We superimpose it to the unstable invariant manifold of the planar Lyapunov 
orbit of the same energy and we observe that they describe approximately the same loci in the
$(x,y)$-projection, the invariant manifold of the invariant torus being located inside the invariant 
manifold of the planar periodic orbit.

\begin{figure}[!ht]
\begin{center}
\includegraphics[scale=0.3,angle=-90]{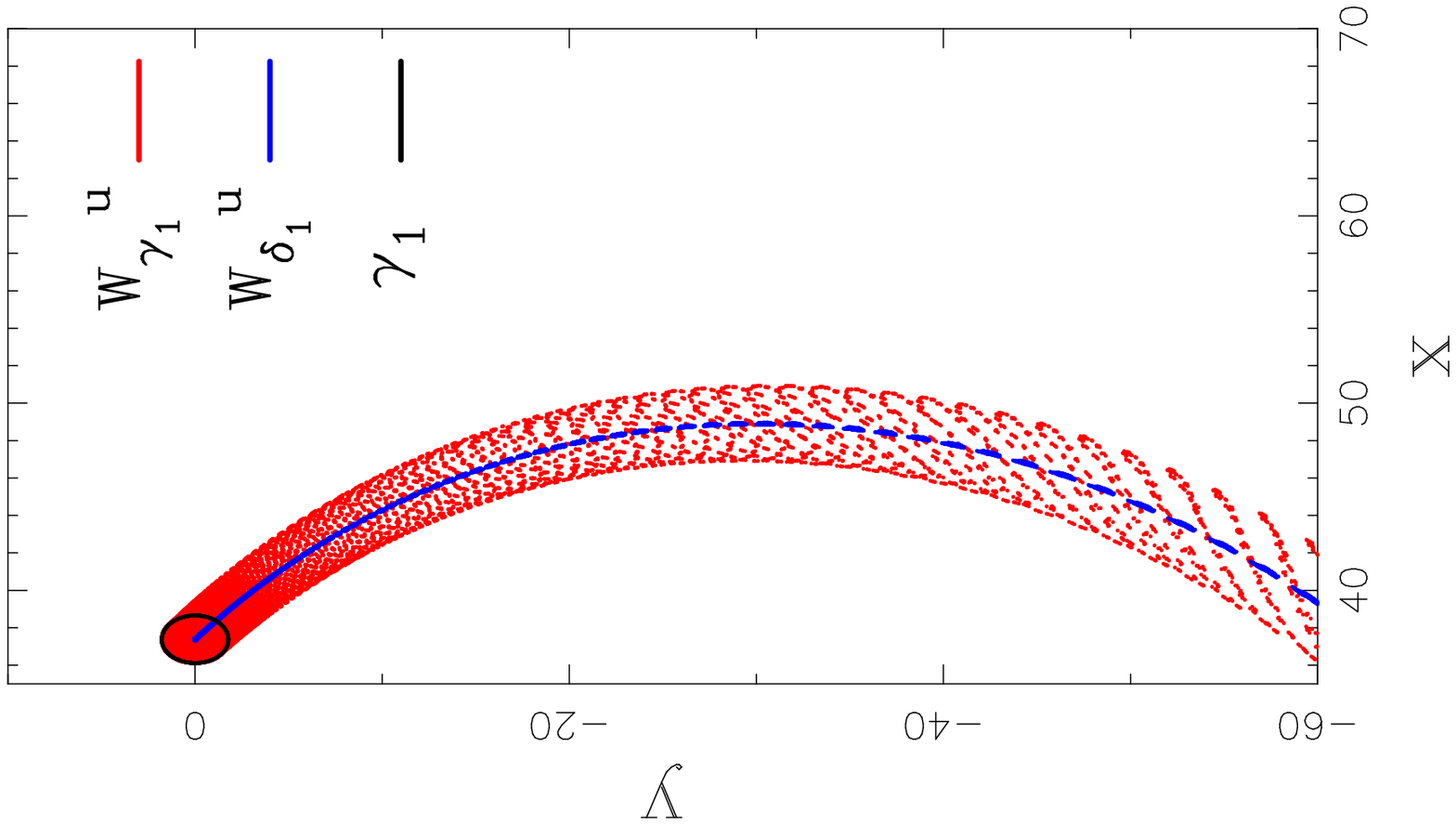} 
\includegraphics[scale=0.3,angle=-90]{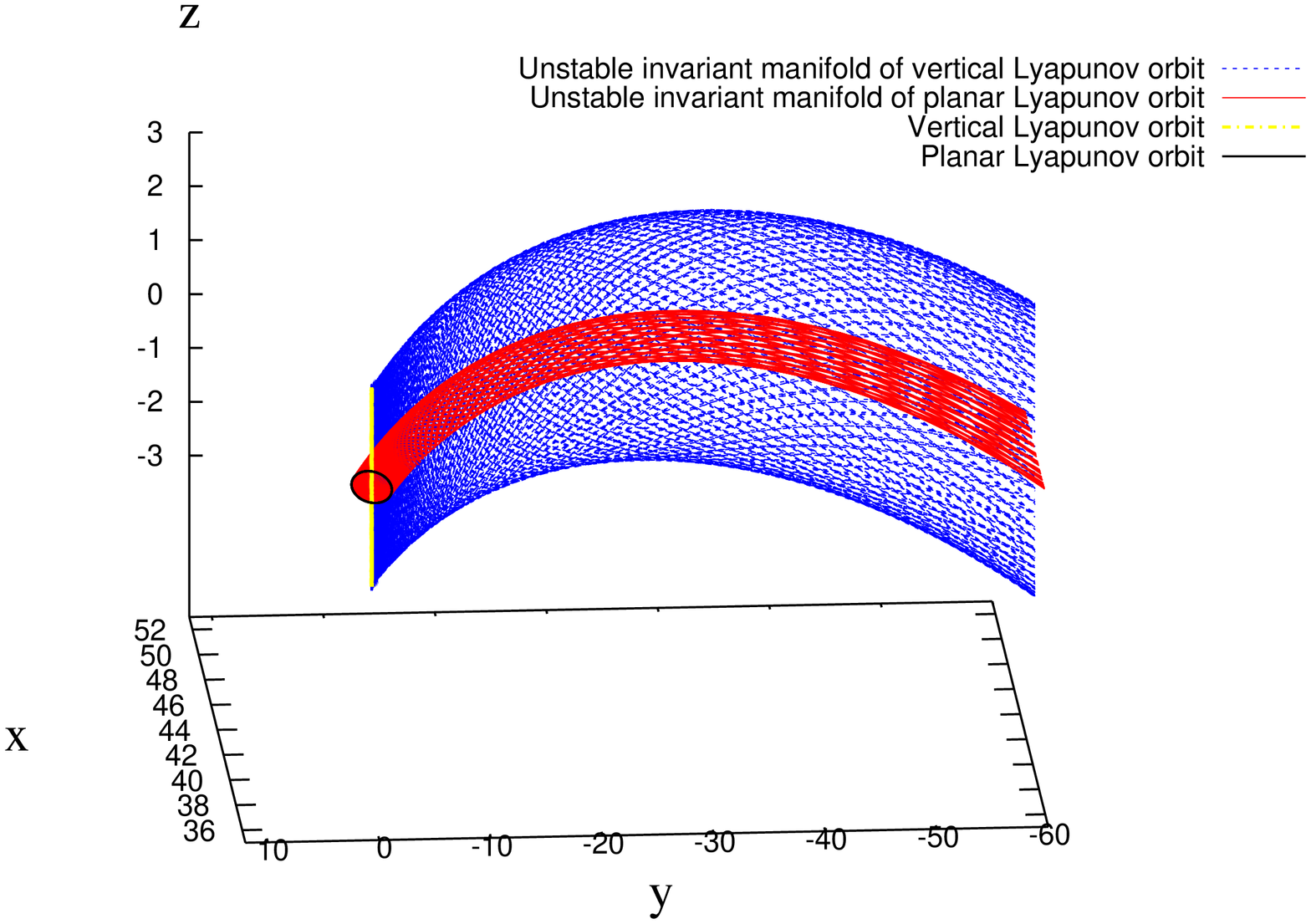}
\end{center}
\caption{Two different views of the invariant manifolds of the planar Lyapunov periodic orbit, 
$W^u_{\gamma_1}$ (in red), and of the vertical Lyapunov periodic orbit, $W^u_{\delta_1}$ (in blue), 
of the logarithmic model for the energy level of $E_J=130155.178$. {\sl Left panel}: $(x,y)$ 
plane. {\sl Right panel:} 3D view (not to scale). The black solid line and the yellow dot-dashed 
line show the planar and the vertical Lyapunov orbits around $L_1$, $\gamma_1$ and $\delta_1$, 
respectively.}
\label{fig:ver-pla}
\end{figure}

In Fig.~\ref{fig:ver-pla}, we plot the unstable invariant manifolds of both 
the planar, $\gamma_1$, and vertical, $\delta_1$, Lyapunov orbits of the same energy 
level ($E_J=130155.178$). The initial conditions used to compute the unstable invariant manifold 
of the vertical Lyapunov orbit have the form $(q_1,p_1,q_2,p_2,q_3,p_3)=(\pm \epsilon,0,0,0,q_3,p_3)$ 
and the initial conditions for the stable invariant manifold are $(q_1,p_1,q_2,p_2,q_3,p_3)$ 
$=(0,\pm \epsilon,0,0,q_3,p_3)$, with $\epsilon=10^{-5}$. In the left panel, we plot the invariant 
manifolds on the projected galactic plane. We observe that they both describe essentially the same 
loci, the invariant manifold of the vertical Lyapunov orbit being located inside the invariant 
manifold of the planar Lyapunov orbit in the configuration space. In the right panel, we show the 
three dimensional view of the same manifolds. We observe that the invariant manifolds of the planar 
Lyapunov orbits drive the motion dominating over the invariant manifolds of the vertical Lyapunov 
orbits. This is due to the fact that the unstable component lies within the plane, since the two real 
eigenvalues are associated with the in-plane motion. We can, therefore, conclude that the motion in the 
neighbourhood of $L_1$ and $L_2$ is basically dominated by the invariant manifolds of planar periodic 
orbits.

\subsection{The role of invariant manifolds in the transfer of matter}
\label{sec:trans}
As we have seen, the motion in a neighbourhood of a hyperbolic point is mainly driven by the 
invariant manifolds of planar periodic orbits, so in this section we restrict ourselves to the
${z=0}$ plane. The next step is to study the global behaviour of the invariant manifolds associated
with planar periodic orbits and determine their role in global structures. In particular, we are 
interested in knowing whether the invariant manifolds can drive particles from the neighbourhood of 
$L_1$ to the neighbourhood of $L_2$ or even connect the neighbourhood of $L_1$ with itself using 
long-time trajectories. For clarity, as we did in previous sections, we plot the Poincar\'e maps 
and the orbits in the initial frame using barycentric coordinates.

\begin{figure}[!ht]
\begin{center}
\includegraphics[scale=0.4,angle=-90]{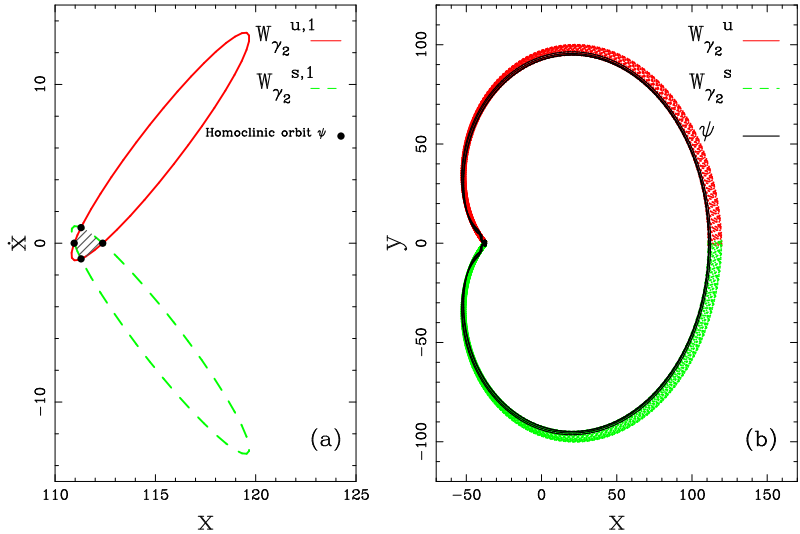}
\caption{Homoclinic connections for a model with $\Omega=5\,\rm{km}\,\rm{s}^{-1}\,\rm{kpc}^{-1}$, 
$v=200\,\rm{km}\,\rm{s}^{-1}$, $R_0=14.14\,\rm{kpc}$, $\p=0.7$ and $\q=0.65$ (Not to scale). 
{\bf (a)} We plot the closed curves $W^{u,1}_{\gamma_2}$ (red solid line) and $W^{s,1}_{\gamma_2}$ 
(green dashed line) on the surface of section $\cal{S}$ defined by $y=0$ with $x>0$. The 
black dots correspond to homoclinic orbits. {\bf (b)} We plot the same invariant manifolds as 
in panel (a), now in the configuration space $(x,y)$. The black curves correspond to the homoclinic 
orbits.}
\label{fig:homo}
\end{center}
\end{figure}

A particle will be transferred from the vicinity of one Lagrangian point to the vicinity of the
symmetric point if it is ``trapped'' first by $W^u_{\gamma_1}$ and then by $W^s_{\gamma_2}$. To study 
these type of transitions, we use Poincar\'e surfaces of section, that is, we draw the crossings of 
trajectories through a particular plane or surface in phase space. Depending on the purposes of 
our study, some surfaces will be more suitable than others, but the methodology is the same. 
Let us take as an example the surface of section $\cal{S}$ defined by $y=0$ with $x>0$; that is,
we consider the orbits when they cut the plane $y=0$ having a positive value
for the $x$ coordinate. Let us consider this surface of section $\cal{S}$ for the stable
and unstable invariant manifolds of a Lyapunov orbit around $L_2$ (located in the $x<0$ side).
Taking initial conditions on the manifold and integrating $W^u_{\gamma_2}$ forward
in time (resp. $W^s_{\gamma_2}$ backwards in time) until the first encounter with
$\cal{S}$, we obtain the simple closed curves $W^{u,1}_{\gamma_2}$ (resp.
$W^{s,1}_{\gamma_2}$) that can be seen in Fig.~\ref{fig:homo}a. Although the simple 
closed curves $W^{u,1}_{\gamma_2}$ and $W^{s,1}_{\gamma_2}$
are obtained as the natural result of intersecting the manifold tubes with a
plane, it should be mentioned that further crossings (i.e. $W^{u,k}_{\gamma_2}$
and $W^{s,k}_{\gamma_2}$ with $k>1$) may well not have this simple structure,
as can be seen in the \rtbpe example provided by Gidea and Masdemont \cite{gid07}.

In the selected example, $W^{u,1}_{\gamma_2}$ and $W^{s,1}_{\gamma_2}$ are represented 
in $(x,\dot{x})$ coordinates. It is important to note that a pair $(x,\dot{x})$ in 
$\cal{S}$ defines an orbit in a unique way, since $y=0$ and $\dot{y}$ is obtained from 
the energy level under study and in the sense of crossing $\cal{S}$. By definition of 
invariant manifold, a point in $W^{u,1}_{\gamma_2} \cap W^{s,1}_{\gamma_2}$ (black dots in 
Fig.~\ref{fig:homo}a) represents a trajectory asymptotic to the Lyapunov orbit $\gamma$ around 
$L_2$ both forward and backward in time and this trajectory is called a homoclinic orbit. In general, 
homoclinic orbits correspond to asymptotic trajectories, $\psi$, such that 
$\psi\in W^u_{\gamma_i}\cap W^s_{\gamma_i},\,i=1,2$. Thus, a homoclinic orbit departs 
asymptotically from the unstable Lyapunov periodic orbit $\gamma$ around $L_i$ and returns 
asymptotically to it, as can be seen in Fig.~\ref{fig:homo}b.

\begin{figure}[!ht]
\begin{center}
\includegraphics[scale=0.4,angle=-90]{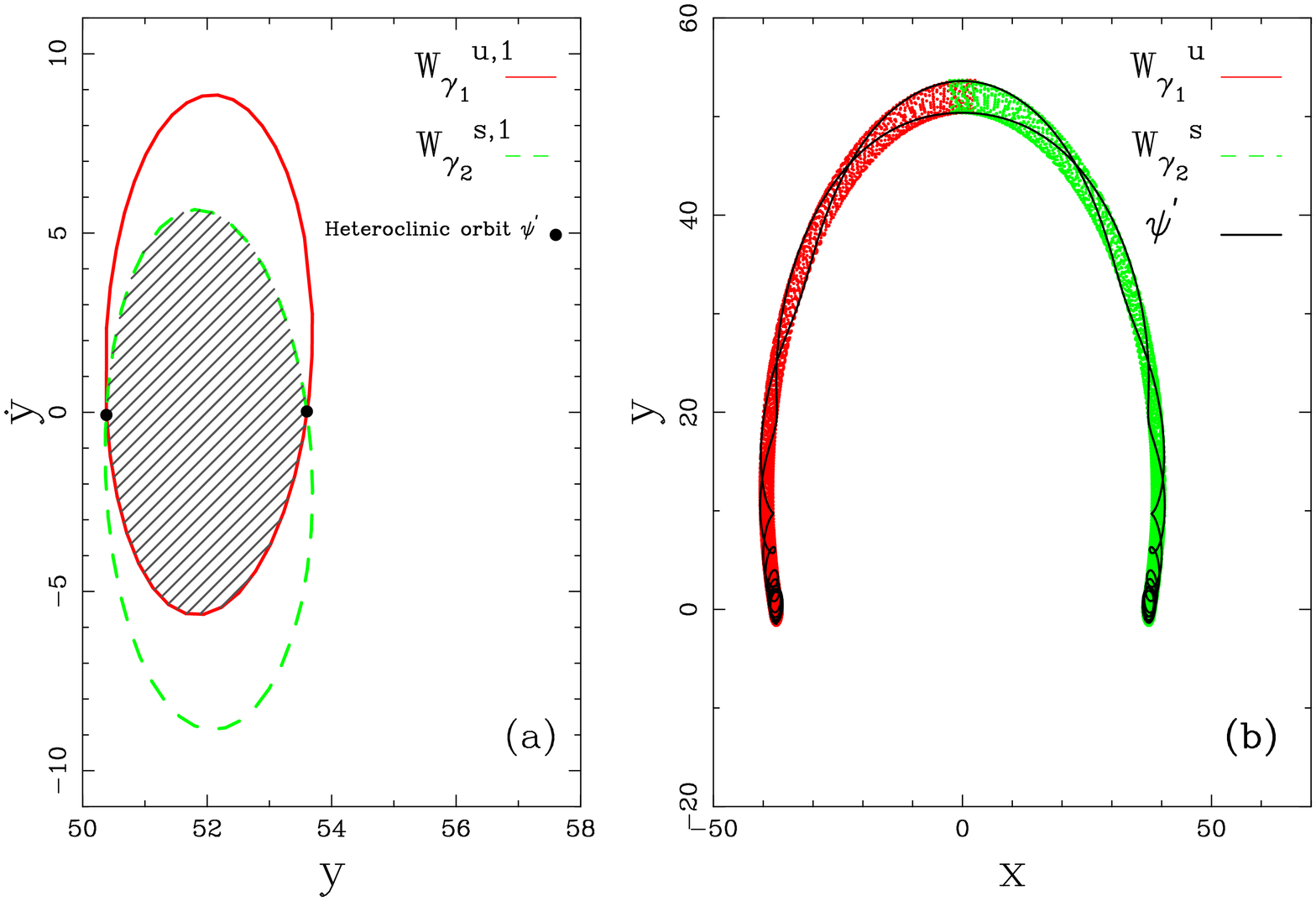}
\caption{Heteroclinic connections for a model with $\Omega=5\,\rm{km}\,\rm{s}^{-1}\,\rm{kpc}^{-1}$,
$v=200\,\rm{km}\,\rm{s}^{-1}$, $R_0=14.14\,\rm{kpc}$, $\p=0.95$ and $\q=0.85$ (Not to scale). 
{\bf (a)} We plot the closed curves $W^{u,1}_{\gamma_1}$ (red solid line) and $W^{s,1}_{\gamma_2}$ 
(green dashed line) on the surface of section $\cal{S}^{\prime}$ defined by $x=0$ with $y>0$. 
The black dots correspond to heteroclinic orbits. {\bf (b)} We plot the same invariant manifolds 
as in panel (a), now in the configuration space $(x,y)$. The black curves correspond to the 
heteroclinic orbits.}
\label{fig:hetero}
\end{center}
\end{figure}

Heteroclinic orbits, on the other hand, are defined as asymptotic trajectories, $\psi^\prime$, 
such that $\psi^\prime\in W^u_{\gamma_i}\cap W^s_{\gamma_j},\, i\ne j,\,i,j=1,2$. Thus,
a heteroclinic orbit departs asymptotically from the periodic orbit $\gamma$ around $L_i$ 
and asymptotically approaches the corresponding Lyapunov periodic orbit of the same 
energy around the Lagrangian point at the opposite end of the bar region $L_j$, $i\ne j$; a
suitable surface of section, $\cal{S}^{\prime}$, for this computation can be the plane $x=0$ with
$y>0$. We hereafter refer to interior region as the elliptical-like shape defined by the iso-effective 
potential curves encircling the centre and roughly passing through $L_1$ and $L_2$. In 
Fig.~\ref{fig:hetero}a, we plot the closed curves $W^{u,1}_{\gamma_1}$ and $W^{s,1}_{\gamma_2}$ on 
the surface $\cal{S}^{\prime}$. Note that there are two intersection points, therefore this model 
presents two heteroclinic orbits that will connect asymptotically one side of the interior region 
with the opposite, as can be seen in Fig.~\ref{fig:hetero}b.

\begin{figure}[!ht]
\begin{center}
\includegraphics[scale=0.9]{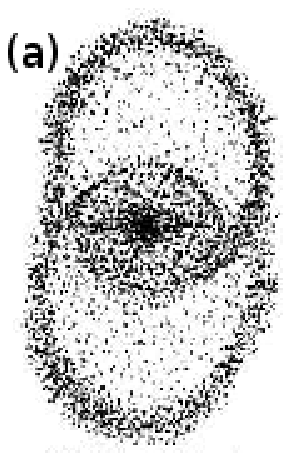}
\includegraphics[scale=1.]{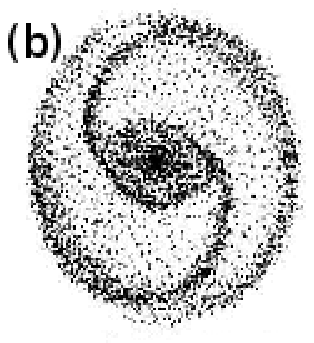}
\includegraphics[scale=1.]{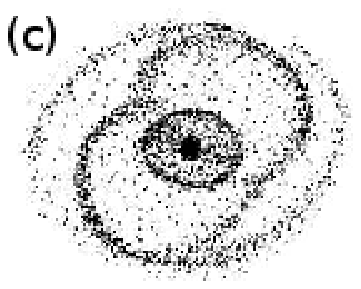}
\includegraphics[scale=0.5]{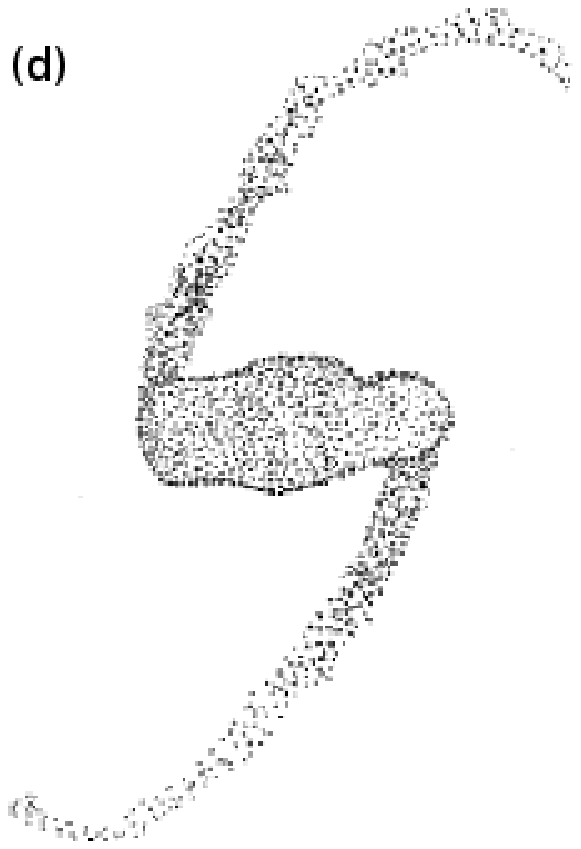}
\caption{{\bf (a), (b)} and {\bf (c)} schematic prototypes of outer rings from Buta \& Crocker \cite{but91}
reproduced by permission of the AAS. {\bf (d)} schematic prototype of a spiral arm.}
\label{fig:schema}
\end{center}
\end{figure}

The existence of these type of connections determines the way matter is transferred within
the galaxy. By definition, if there are heteroclinic connections, they create a tube, 
where particles can be trapped and travel from one side of the interior region to the opposite. The
intersection of the tube with the section $\cal{S}^{\prime}$ can be seen in Fig.~\ref{fig:hetero}a as 
the hatched area obtained from intersecting the curves $W^{u,1}_{\gamma_1}$ and $W^{s,1}_{\gamma_2}$.
The global morphology in such models is reminiscent of that of $R_1$ rings \cite{but91,rom06,rom07}. 
$R_1$ rings are outer rings in barred galaxies whose semi-major axis is perpendicular to the one of 
the bar. They have a characteristic $\theta$ or $8$-shape (see Fig.~\ref{fig:schema}a). If the model presents 
homoclinic connections, particles can be trapped in tubes travelling from one end of the interior 
region to itself. The global morphology in this case is reminiscent of that of $R_1R_2$ rings 
\cite{but91,rom07} (see Fig.~\ref{fig:schema}c). $R_1R_2$ ringed galaxies present two outer rings, 
namely the $R_1$ ring previously mentioned and an $R_2$ ring, whose semi-major axis is parallel to 
the bar semi-major axis (see Fig.~\ref{fig:schema}b). Note that if neither heteroclinic nor homoclinic 
connections exist, particles are not able to return to the vicinity of the interior region in an early 
stage. Thus the particles follow an escaping trajectory forming spiral arms \cite{rom07} (see
Fig.\ref{fig:schema}d).

\begin{figure}[!ht]
\begin{center}
\includegraphics[scale=0.8,angle=-90]{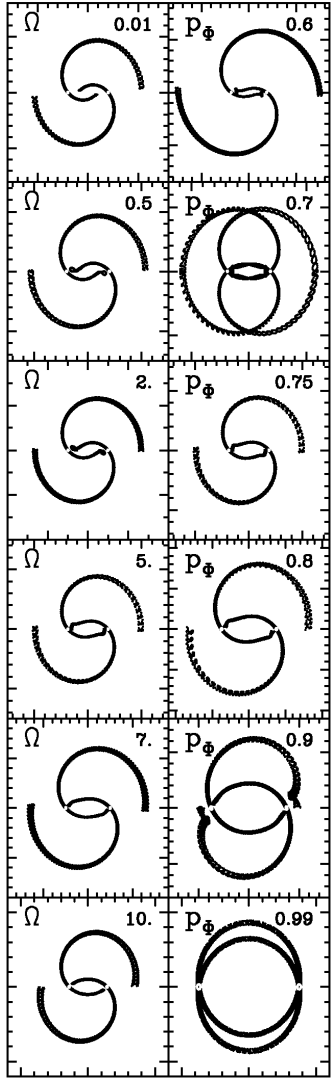}
\caption{Effect of the pattern speed $\Omega$ (left panels) and of the planar shape parameter $\p$ 
(right panels) on the invariant manifolds.}
\label{fig:shape}
\end{center}
\end{figure}

\begin{figure}[!ht]
\begin{center}
\includegraphics[scale=0.3,angle=-90]{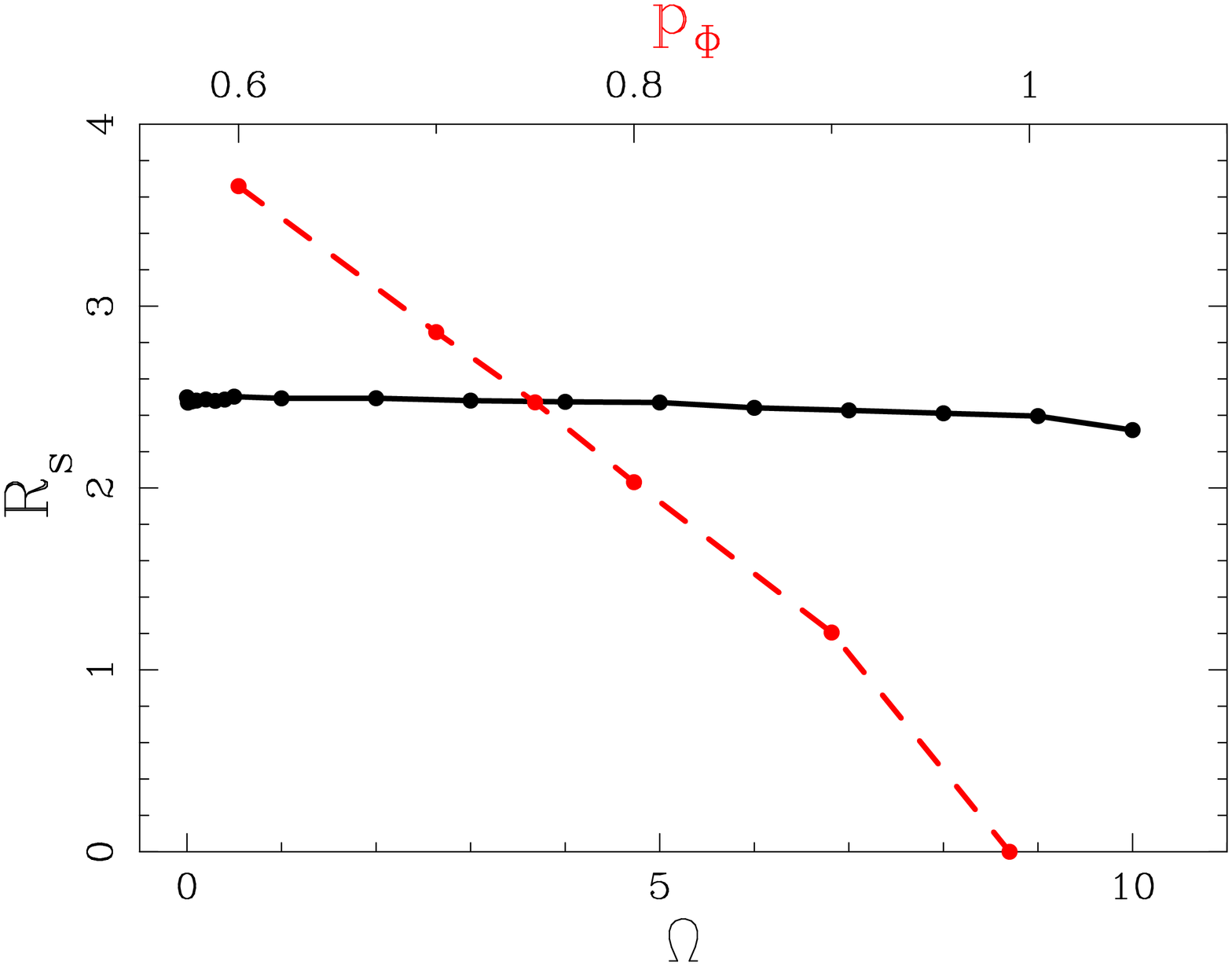}
\caption{Effect of the pattern speed (black solid line and labels in the bottom) and of the 
shape parameter $\p$ (red dashed line and labels in the top) on the degree of openness of 
the invariant manifolds, $R_s$.}
\label{fig:varomega2}
\end{center}
\end{figure}

In Sect.~\ref{sec:invman} we show that motion is mainly driven by the invariant manifolds
of the planar Lyapunov orbits. Here we follow the study in Sect.~\ref{sec:var} and we
compute the invariant manifolds of families of models where we vary the pattern speed, 
$\Omega$, and the planar axial ratio, $\p$, i.e. the two dynamical parameters that have an 
effect on the potential in the galactic plane. The rest of the parameters are fixed as in
``model 1''. For each model, we compute the family of planar Lyapunov orbits in a range 
of energies for which they are unstable and the bottleneck through the zero velocity curve is opened. 
We compute the invariant manifolds associated with the planar periodic orbits and we study 
the existence of possible heteroclinic and homoclinic connections. In each panel of 
Fig.~\ref{fig:shape}, we plot the corresponding global morphology. In the left column we 
increase from top to bottom the value of the pattern speed. In all cases the global shape 
of the invariant manifolds is two spiral arms and the degree of openness is essentially 
independent of the value of the pattern speed. This can be seen in Fig.~\ref{fig:varomega2}, where 
we plot a measure of the openness of the spiral, namely the ratio, $R_s$, between the 
$x$-coordinate of the outer branch of the unstable invariant manifold at the first cut 
with the $y=0$ axis and the absolute value of $x$-coordinate of the equilibrium point. In the 
right column of Fig.~\ref{fig:shape}, we plot the invariant manifolds for models with 
different values of $\p$. When the system is strongly non-axisymmetric, no homoclinic or
heteroclinic connections are present and the global structure is that of two spiral arms. 
When $\p=0.7$, the model presents homoclinic connections with both ends of the 
interior region. As previously mentioned, the global morphology reminds that of $R_1R_2$ rings.
As $\p$ approaches to $1$, that is, the system becomes prolate, the invariant manifolds tend 
to close, heteroclinic connections are present and the structure forms an $R_1$ ring. 
We stress the similarity of these patterns with the schematic prototypes of outer rings 
described by Buta \& Crocker \cite{but91}. In Fig.~\ref{fig:varomega2}, we also plot the 
degree of openness as the shape parameter $\p$ varies. As is illustrated in the right column of 
Fig.~\ref{fig:shape}, the value of the ratio decreases as the value of $\p$ approaches to $1$.

\section{Conclusions}
\label{sec:conc}
In this paper we use a suitable galactic potential (the logarithmic potential) to 
perform a semi-analytic study of the neighbourhood of the unstable equilibrium points
$L_1$ and $L_2$. For a detailed study of their neighbourhood, we use a partial normal
form scheme and we find that the main objects are the planar and vertical families of
Lyapunov orbits and invariant tori. We compute the invariant manifolds associated with both periodic 
orbits and quasi-periodic orbits using the reduced Hamiltonian and we find that the motion around 
$L_1$ and $L_2$ is mainly driven by the invariant manifolds of the planar Lyapunov orbits. However, 
we are also interested in the global structure of the galaxy. Thus, we study the possible homoclinic 
and heteroclinic connections between the planar periodic orbits. For such a purpose, we use suitable 
Poincar\'e surfaces of section. We note that this approach has successfully been used in
celestial mechanics and in this paper we apply it to a galactic dynamics problem, namely 
the formation of spiral arms and rings in barred galaxies, in comparison to other theories given
so far.

Here we are interested in the hyperbolic behaviour of $L_1$ and $L_2$ and, particularly,
in determining the role the invariant manifolds play in the transfer of matter for 
energy levels where the zero velocity curves are open (i.e. a range of energies somewhat larger
than the energy of $L_1$ and $L_2$) and a bottleneck appears around $L_1$ and $L_2$. For
such purpose, we also study the influence of the main model parameters, namely the
pattern speed, $\Omega$, and the shape parameters, $\p$ and $\q$. The logarithmic models are 
suitable for describing triaxial systems such as haloes, bars, bulges in disc galaxies or elliptical 
galaxies. The structures we have constructed using the invariant manifolds, however, are not 
globally dependent on the model characteristics \cite{rom06,rom07}. This implies, for example, that if 
a system had some degree of rotation, these kind of structures should be present. Elliptical galaxies
are triaxial systems that do not present any external feature, i.e. they are an ellipsoidal
distribution of matter with different degrees of ellipticity and they present neither
spiral arms nor rings. Observations show that elliptical galaxies barely rotate or do not rotate
at all as a figure \cite{bin87,sch82}. Our results are in agreement with this statement. We have shown 
that models that rotate slowly or do not rotate cancel the hyperbolic behaviour of the equilibrium 
points and, thus, no transfer or escape of matter is possible. On the other hand, bars in disc galaxies are 
non-axisymmetric components usually characterised in the literature by elliptical distributions of 
density \cite{pfe84} although both observations and simulations show they might have more rectangular 
endings \cite{ath90,ath02}. It is well-known that bars rotate and observations show that spiral arms 
or rings emanate from the ends of the bar. This characteristic is also consistent with our results which 
show that if the system rotates at a given angular velocity, the hyperbolic equilibrium points are 
present, so that the invariant manifolds drive the motion, and therefore, set the global morphology
to the galaxy. We have seen that, depending on the rotation velocity and the shape of the bar, the
morphology will be that of a barred spiral galaxy or that of a barred ringed galaxy.

\section*{Acknowledgements}
This work partially supported by the Spanish grants MCyT-FEDER MTM2006-00478 and AYA2007-60366 and 
the French grant ANR-06-BLAN-0172. MRG acknowledges her ``Becario MAE-AECI'' and the Marie Curie
Research Training Network Astronet. 

\bibliographystyle{plainnat}

\end{document}